\documentclass[12pt,reqno]{amsart}

\usepackage{amsmath,amsfonts,amsthm,amssymb, latexsym}
\usepackage[T1]{fontenc}
\usepackage[latin1]{inputenc}

\newcommand{\bp}{\mathbf{P}}
\newcommand{\be}{\mathbf{E}}
\newcommand{\bd}{\mathbf{\Gamma}^{\mathbf{2}}}
\newcommand{\bde}{\mathbf{\Gamma}^{\mathbf{2},\ep}}
\newcommand{\bdet}{\mathbf{\Gamma}^{\mathbf{2},\eta}}
\newcommand{\bx}{{\bf X}}
\newcommand{\der}{\delta}
\newcommand{\dom}{\mbox{Dom}}

\newcommand{\ha}{\hat a}

\newcommand{\hr}{\hat r}
\newcommand{\hsi}{\hat \sigma}
\newcommand{\hz}{\hat z}
\newcommand{\ka}{\kappa}
\newcommand{\im}{\imath}

\newcommand{\ibst}{\int_{[\bar s,\bar t]}}
\newcommand{\ist}{\int_{[s,t]}}
\newcommand{\norm}[1]{\lVert #1\rVert}

\newcommand{\ott}{[0,T]}

\newcommand{\pp}{\Pi^{+}}
\newcommand{\sgn}{{\rm sgn}}
\newcommand{\xd}{{\bf X^{2}}}
\newcommand{\xdst}{{\bf X}_{ts}^{\bf 2}}

\newcommand{\half}{{\frac{1}{2}}}

\DeclareMathOperator{\id}{\text{Id}}

\newcommand{\ca}{{\mathcal A}}
\newcommand{\cb}{{\mathcal B}}
\newcommand{\cac}{{\mathcal C}}

\newcommand{\cf}{{\mathcal F}}

\newcommand{\cj}{{\mathcal J}}

\newcommand{\cn}{{\mathcal N}}

\newcommand{\cq}{{\mathcal Q}}
\newcommand{\crr}{{\mathcal R}}

\newcommand{\cu}{{\mathcal U}}
\newcommand{\cv}{{\mathcal V}}

\newcommand{\cZ}{{\mathcal Z}}
\newcommand{\cz}{{\mathcal Z}}

\newcommand{\al}{\alpha}
\newcommand{\ga}{\gamma}
\newcommand{\gga}{\Gamma}
\newcommand{\ep}{\varepsilon}
\newcommand{\eps}{\epsilon}

\newcommand{\om}{\omega}
\newcommand{\si}{\sigma}
\newcommand{\vp}{\varphi}
\newcommand{\ze}{\zeta}

\newcommand{\laa}{\Lambda}
\newcommand{\oom}{\Omega}

\newcommand{\C}{{\mathbb C}}

\newcommand{\N}{{\mathbb N}}

\newcommand{\R}{{\mathbb R}}

\newcommand{\lcl}{\left\{}
\newcommand{\rcl}{\right\}}
\newcommand{\lp}{\left(}
\newcommand{\rp}{\right)}
\newcommand{\lc}{\left[}
\newcommand{\rc}{\right]}
\newcommand{\lln}{\left|}
\newcommand{\rrn}{\right|}

\newcommand{\bean}{\begin{eqnarray*}}
\newcommand{\eean}{\end{eqnarray*}}
\newcommand{\ben}{\begin{enumerate}}
\newcommand{\een}{\end{enumerate}}
\newcommand{\beq}{\begin{equation}}
\newcommand{\eeq}{\end{equation}}

\newtheorem{theorem}{Theorem}[section]

\newtheorem{definition}[theorem]{Definition}

\newtheorem{hypothesis}{Hypothesis}
\newtheorem{lemma}[theorem]{Lemma}
\newtheorem{notation}[theorem]{Notation}

\newtheorem{proposition}[theorem]{Proposition}

\theoremstyle{remark}
\newtheorem{remark}[theorem]{Remark}

\begin{document}

\title[analytic fractional Brownian motion]{The rough path associated to the multidimensional analytic fbm with any Hurst parameter}
\author[S. Tindel \and J. Unterberger]{Samy Tindel \and Jérémie Unterberger}

\address{
{\it Samy Tindel and Jérémie Unterberger:}
{\rm Institut {\'E}lie Cartan Nancy, Nancy-Universit\'e, B.P. 239,
54506 Vand{\oe}uvre-l{\`e}s-Nancy Cedex, France}.
{\it Email: }{\tt tindel@iecn.u-nancy.fr}, {\tt unterber@iecn.u-nancy.fr}
}



\date{\today}

\begin{abstract}
In this paper, we consider a complex-valued $d$-dimensional fractional Brownian motion defined on the closure of the complex upper half-plane, called {\it analytic fractional Brownian motion} and denoted by $\gga$. This process has been introduced in \cite{Un}, and both its real and imaginary parts, restricted on the real axis, are usual fractional Brownian motions. The current note is devoted to prove that a rough path based on $\gga$ can be constructed for any value of the Hurst parameter in $(0,1/2)$. This allows in particular to solve differential equations driven by $\gga$ in a neighborhood of 0 of the complex upper half-plane, thanks to a variant of the usual rough path theory due to Gubinelli \cite{Gu}.
\end{abstract}
\maketitle


\section{Introduction}

The (two-sided) fractional Brownian motion $t\to B_t$, $t\in\R$  (fBm for short) with Hurst exponent $\alpha$, $\alpha\in(0,1)$, defined as the centered Gaussian process with covariance
\begin{equation}\label{eq:1} 
\be[B_s B_t]=\half (|s|^{2\alpha}+|t|^{2\alpha}-|t-s|^{2\alpha}), \end{equation}
is a natural generalization in the class of Gaussian processes of
the usual Brownian motion, in the sense that it exhibits two fundamental properties shared with Brownian motion, namely,
it has stationary increments, viz. $\be[(B_t-B_s)(B_u-B_v)]=\be[(B_{t+a}-B_{s+a})(B_{u+a}-B_{v+a})]$ for
every $a,s,t,u,v\in\R$, and it is self-similar, viz. 
\begin{equation} \label{eq:2} 
\forall \lambda>0, \quad (B_{\lambda t}, t\in\R) \overset{(law)}{=} (\lambda^{\alpha} B_t,
t\in\R).
\end{equation} 
One may also define  a $d$-dimensional 
vector Gaussian process (called: {\it $d$-dimensional fractional Brownian motion}) by setting $B_t:={B}_t=(B_t(1),\ldots,B_t(d))$, where $(B_t(i),t\in\R)_{i=1,\ldots,d}$ are $d$ independent (scalar) fractional Brownian motions.
Its theoretical interest lies in particular in the fact that it is (up to normalization) the only Gaussian process satisfying the two properties (\ref{eq:1}) and (\ref{eq:2}).
Furthermore, a standard application of Kolmogorov's theorem shows that fBm has a version with $(\alpha-\eps)$-H\"older paths for every $\eps>0$. This makes this process amenable to models where a Gaussian process with Hölder continuity exponent different from $1/2$ is needed, and we refer for instance to \cite{CLCHH,KS,OTHP} for some applications to biophysics.

\smallskip

Consequently, there has been a widespread interest during the past ten years in constructing a stochastic integration theory
with respect to fBm and solving stochastic differential equations driven by fBm. The multi-dimensional case is
very different from the one-dimensional case. When one tries to integrate for instance a stochastic differential
equation driven by a two-dimensional fBm ${B}=({B}(1),{B}(2))$ by using any kind of Picard iteration scheme, one
encounters very soon the problem of defining the L\'evy area of ${B}$ which is the antisymmetric part
of ${\ca}_{ts}:=\int_s^t dB_{t_1}(1) \int_s^{t_1} dB_{t_2}(2)$. This is the simplest occurrence
of iterated integrals ${B}^k_{ts}(i_1,\ldots,i_k):=\int_s^t dB_{t_1}(i_1)\ldots \int_s^{t_{k-1}} dB_{t_k}(i_k)$
which lie at the heart of the rough path method due to T. Lyons. 

\smallskip

Let us describe briefly this method, rephrased in the setting of \cite{Gu} which is going to be used in the sequel of the paper: assume
${X}=({X}(1),\ldots,{X}(d))$ is some non-smooth $\alpha$-H\"older $d$-dimensional path. Integrals such as
$\int f_1({X}(t))d{X}_1(t)+\ldots+f_d({X}(t))d{X}_d(t)$ do not make sense a priori because
${X}$ is not differentiable (Young's integral works for $\alpha>\half$ but not beyond). In order to define
the integration of a differential form along ${X}$, it is enough to define a truncated {\it multiplicative functional}
$({X}^{1},\ldots,{X}^{\lfloor 1/\al \rfloor})$ where ${X}^{1}_{ts}=X_t-X_s$ and
each
${X}^k=({X}^k(i_1,\ldots,i_k))_{1\le i_1,\ldots,i_k\le d}$ -- a matrix of (increments of)
continuous paths --
is a {\it substitute} for the iterated integrals formally given as ${X}^k_{ts}:=\int_s^t d{X}_{t_1}(i_1)\int_s^{t_1} d{X}_{t_2}(i_2)
\ldots \int_{s}^{t_{k-1}} d{X}_{t_k}(i_k)$, with the following two properties:

\begin{itemize}
\item[(i)] Each component of ${X}^k$ is $k\ka$-Hölder continuous for any $\ka<\al$.
\item[(ii)] {\it Multiplicativity}: letting $(\der X^k)_{tus}:=X^k_{ts}-X^k_{tu}-X^k_{us}$, one requires
\begin{equation}
(\der X^k)_{tus}(i_1,\ldots,i_k)=\sum_{k_1+k_2=k} 
{X}_{tu}^{k_1}(i_1,\ldots,i_{k_1}) {X}_{us}^{k_2}(i_{k_1+1},\ldots,i_k). \label{eq:*}
\end{equation}
\end{itemize}
Once these functionals are defined, the theory described in \cite{FV,Gu,LQ-bk} can be seen as a procedure which allows to define out of these data iterated integrals of any order and
to solve differential equations driven by ${X}$.

\smallskip

With these preliminary considerations in mind, it is easily conceived that the fundamental problem in order to apply the general theory is to give a suitable definition of the functionals $X^k$. For any smooth path, $X^k$ can be defined as a Riemann multiple integral. The multiplicative and Hölder continuity properties are then trivially satisfied, as can be checked by direct computation.
 So the most
natural way to construct such a multiplicative functional is to start from some smooth approximation
${X}^{\eta}$, $\eta\searrow 0$ of ${X}$ such that each iterated
 integral ${X}^{k,\eta}(i_1,\ldots,i_k)$, $k\le \lfloor 1/\al \rfloor$ converges in the
 $k \ka$ -Hölder norm for any $\ka<\al$. This general scheme has been applied to fBm in a paper by L. Coutin and Z. Qian \cite{CQ}, by means of  standard $n$-dyadic piecewise
 linear approximations $\tilde{B}^{2^{-n}}$ of $B$. In a later paper, one of the authors \cite{Un} tried to tackle the problem by  seeing ${B}$ as the real part of
 the boundary value of 
an analytic process  $\gga$  
living on the upper half-plane $\Pi^+=\{z\in\C\ | \Im z>0\}$. The time-derivative of this
centered Gaussian process has the following hermitian positive-definite covariance kernel:
\begin{equation}\label{eq:3} 
\be \lc {\Gamma}'(z)\overline{{\Gamma}'(w)} \rc= K'(z,\bar{w})=
\frac{\alpha(1-2\alpha)}{2\cos\pi\alpha}
(-\imath(z-\bar{w}))^{2\alpha-2},\quad z,w\in\Pi^+\end{equation}
where  $z^{2\alpha-2}:=e^{(2\alpha-2)\ln z}$ (with the usual determination of the logarithm)
is defined and analytic on the cut plane $\C\setminus\R_-$. Also, by construction, 
$\be {\Gamma}'(z){\Gamma}'(w)\equiv 0$ identically. It is essential to understand that $K'$ is a 
multivalued function on $\C\times\C\setminus\{(z,\bar{z})\ |\  z\in\C\}$; on the other hand, for $z,w\in\Pi^+$ we have 
$\Re(-\imath(z-\bar{w}))>0$, so the kernel $K'$ is well-defined.
Then $\Re B_t^{\eta}:=\Re {\Gamma}_{t+\imath\eta}$ is a good approximation of fBm, namely, $\Re B^{\eta}$
 converges a.s. in
the $\ka$-Hölder distance to a process $\Re B$ with the same law as fBm for any $\kappa<\alpha$.

\smallskip

Both approximation schemes introduced in \cite{CQ,Un} lead to the same semi-quantitative result, namely:
\begin{itemize}
\item
When $\alpha>1/4$, the L\'evy area and volume (in other words, the truncated
 multiplicative functional truncated
to order $3$) converge a.s. in the appropriate variation norm. The heart of the proof
 lies in the study of the L\'evy
area $\tilde{\ca}_{ts}^{2^{-n}}$, resp. ${\ca}_{ts}^{\eta}$ of the smooth
 approximation; one may prove in
particular that $\be[(\tilde{\ca}_{ts}^{2^{-n}})^2]$ and 
$\be[({\ca}_{ts}^{\eta})^2]$
converge to the same limit when $2^{-n}$ or $\eta$ go to 0;
\item
When $\alpha<1/4$, $\be[(\tilde{\ca}_{ts}^{2^{-n}})^2]$ and
 $\be[({\ca}_{ts}^{\eta})^2]$
diverge resp. as $n^{(1-4\alpha)}$ and $\eta^{-(1-4\alpha)}$. Hence the methods alluded to above fail.
\end{itemize}
The latter result is of course unsatisfactory, and constitutes by no means a proof that no coherent
stochastic integration theory
with respect to fBm may exist when $\alpha<1/4$. However, to the best of our knowledge, there is no explicit example in the literature of a $d$-dimensional (with $d>1$) process with Hölder regularity $\ka<1/4$ allowing the construction of a rough path.

\smallskip

The current article proposes then to make a step in this direction, and the rough path construction we propose will simply be obtained by considering the complex-valued process $\Gamma$ (recall that this process is induced by the covariance function (\ref{eq:3})) for its own sake, instead of $\Re\Gamma$. In particular, for $t\in\R$, the irregular process $B_t:=\Gamma_t$ will be approximated by the complex valued analytic process ${B}_t^{\eta}:={\Gamma}_{t+\imath\eta}$. But in a more general way, $\gga^{\eta}_{z}:=\gga_{z+\imath \eta}$ will stand for an analytic approximation of $\gga$ on the closed upper half-plane $\bar\Pi^+:=\{z\in\C\ |\ \Im z\ge 0\}$.
 An adequate  limiting procedure for $\eta\to 0$
 will allows us to  prove the following main results:

\begin{enumerate}
\item
The iterated integrals ${\Gamma}^{\eta,k}(i_1,\ldots,i_k)$, $k\le \lfloor 1/\al \rfloor$
converge in the $(k\ka)$-Hölder norm for any $\ka<\al$ and {\em any} Hurst index $\alpha>0$. The limiting objects $\Gamma^k$ satisfy our conditions (i) and (ii) above, which yields the construction of a rough path above the process $\Gamma$.
\item
One deduces from this fact that
stochastic differential equations of the type
\begin{equation} dz_t=b(z_t)dt+\sigma(z_t)d{\Gamma}_t,\ \quad z_0=a\in\C^n, \label{eq:SDE} \end{equation}
where $b$ and $\sigma$ are vector-valued, resp. matrix-valued analytic functions on a complex
neighborhood $\Omega$ of $0$, have a strong local solution $z_t$ defined on
$\Omega'\cap \bar{\Pi}^+$ where $\Omega'\subset \Omega$ is another complex neighborhood of $0$.
\end{enumerate}
Let us make now a few comments on these results:

\smallskip

\noindent
{\it (i)} 
An appropriate name for the ${\Gamma}$-process could be {\em analytic fractional Brownian motion}
(analytic fBm or afBm for short). This is the name we shall use throughout the article. Yet the reader should
be warned against two possible misunderstandings:

\noindent
-- ${\Gamma}$ is analytic only in the (open) upper half-plane. When we consider its restriction to
$\R$ (its boundary value on $\R$, one might say) it is merely a continuous process with the same H\"older
continuity as the usual fBm. The fact that ${\Gamma}$ is very irregular on $\R$ makes it interesting
to be able to solve stochastic differential equations driven by ${\Gamma}$, whereas they are almost
trivially solved on $\Pi^+$;

\noindent
--
considering the restriction of ${\Gamma}$ to $\R$, one may be tempted to write ${\Gamma}_t=\Re {\Gamma}_t+\imath
\Im{\Gamma}_t$ and to consider separately the real and the imaginary part. Elementary computations
show that both $\Re {\Gamma}$ and $\Im {\Gamma}$ have the same law as fBm. But ${\Gamma}$ is not merely a
{\em complex} fBm since $\Re {\Gamma}$ and $\Im {\Gamma}$ are {\em not} independent (see section 1).
It is the correlation between $\Re {\Gamma}$ and $\Im {\Gamma}$ that cancels the singularities for
small Hurst indices.

\smallskip

\noindent
{\it (ii)} 
It is of course possible to consider (\ref{eq:SDE}) as a system of
two {\em real} coupled equations on $\R^n$, namely,
\begin{eqnarray}
&& dx_t=b_1(x_t,y_t) dt+\sigma_1(x_t,y_t)d\Re {\Gamma}_t-\sigma_2(x_t,y_t)d\Im {\Gamma}_t, \nonumber\\
&& dy_t=b_2(x_t,y_t) dt+\sigma_1(x_t,y_t)d\Im{\Gamma}_t+\sigma_2(x_t,y_t)d\Re {\Gamma}_t
\end{eqnarray}
where $z_t=x_t+\imath y_t$, $b=b_1+\imath b_2$, $\sigma=\sigma_1+\imath\sigma_2$, $\Re {\Gamma}$ and $\Im {\Gamma}$
are correlated fBm, and the components of $(b_1,b_2)$ and $(\sigma_1,\sigma_2)$ satisfy the Cauchy-Riemann
equations $\partial_x b_1=\partial_y b_2$, $\partial_y b_1=-\partial_x b_2$ and 
$\partial_x \sigma_1=\partial_y \sigma_2$, $\partial_y \sigma_1=-\partial_x \sigma_2$.

\smallskip

\noindent
{\it (iii)} As in \cite{GT,TT}, we have chosen here to solve the differential equation (\ref{eq:SDE}) thanks to a variant of Lyon's rough path method \cite{FV,LQ-bk} called {\em algebraic integration
theory}, introduced in  \cite{Gu}. The ideas are roughly the same in both theories, but the technical apparatus
 is  different and makes some of our technical proofs become significantly shorter (in particular,
one replaces $q$-variation norms by H\"older norms which are easier to estimate), and renders the adaptation of the theory to the analytic setting more transparent. We outline this generalization to the complex plane in Sections \ref{sec:alg-integration}, \ref{sec:rough-2} and \ref{sec:rough-general}.

\smallskip

\noindent
{\it (iv)}
Let us try to explain briefly why the regularized Levy area is divergent for the real fBm and $\al\le 1/4$, while it converges for the analytic one for any $\al>0$. Let us call then ${\ca}_{ts}^{\eta}$ the regularized Levy area for $\Re\Gamma$, and let us compute $\be ({\ca}_{ts}^{\eta})^2$ for $t,s\in\R$: by definition (recall that $\be[ {\Gamma}'(z){\Gamma}'(w)] \equiv 0$
identically)
\begin{eqnarray} \be \lc ({\ca}_{ts}^{\eta})^2\rc &=& 2\be \left( \int_s^t d{\Gamma}_{x_1+\imath\eta}(1) \int_s^{x_1}
d{\Gamma}_{x_2+\imath\eta}(2) \right)\left( \int_s^t d\bar{{\Gamma}}_{y_1+\imath\eta}(1)\int_s^{y_1}
d\bar{{\Gamma}}_{y_2+\imath\eta}(2) \right) \nonumber\\
&+& 2\Re \be  \left( \int_s^t d{\Gamma}_{x_1+\imath\eta}(1) \int_s^{x_1}
d\bar{{\Gamma}}_{x_2+\imath\eta}(2) \right)\left( \int_s^t d\bar{{\Gamma}}_{y_1+\imath\eta}(1)\int_s^{y_1}
d{\Gamma}_{y_2+\imath\eta}(2) \right) \nonumber\\
&=:& {\cv}_1(\eta)+{\cv}_2(\eta).
\end{eqnarray}
The first term in the right-hand side writes
\begin{eqnarray*}
&& {\cv}_1(\eta)\\
&&=C \int_s^t dx_1 \int_s^{x_1} dx_2 \int_s^t dy_1 \int_s^{y_1} dy_2
(-\imath(x_1-y_1)+2\eta)^{2\alpha-2} (-\imath(x_2-y_2)+2\eta)^{2\alpha-2}\\
&& = C' \int_s^t dx_1 \int_s^t dy_1 (-\imath(x_1-y_1)+2\eta)^{2\alpha-2} \\
&& \hspace{4cm}\times\left[(-\imath(x_1-y_1)+2\eta)^{2\alpha}-
(-\imath x_1+2\eta)^{2\alpha}-(\imath y_1+\eta)^{2\alpha} \right],\\
\end{eqnarray*}
while the second term writes
\begin{multline*}
{\cv}_2(\eta)=C' \int_s^t dx_1 \int_s^t dy_1 (-\imath(x_1-y_1)+2\eta)^{2\alpha-2}  \\
\times  \left[(\imath(x_1-y_1)+2\eta)^{2\alpha}-
(\imath x_1+2\eta)^{2\alpha}-(-\imath y_1+\eta)^{2\alpha} \right] .
\end{multline*}
Both integrals look the same {\em except} that ${\cv}_2$ (contrary to ${\cv}_1$) involves
both $-\imath x_1$ and $\imath x_1$, and similarly for $y_1$. This seemingly insignificant difference
is essential, since ${\cv}_1$ can be shown to have a bounded limit when $\eta\to 0$ by using a 
contour deformation in $\Pi^+\times\Pi^+$ which avoids the real axis where singularities live,
while this is impossible for ${\cv}_2$. Namely, $(-\imath(x_1-\bar y_1)+2\eta)^{2\alpha-2}$ is well-defined
if $(x_1,y_1)$ are in the closure of $\Pi^+\times\Pi^+$, while $(\imath(x_1-\bar y_1)+2\eta)^{2\alpha}$
for instance is well-defined on the closure of $\Pi^-\times\Pi^-$, where $\Pi^-$ is the 
{\em lower} half-plane. In fact, explicit computations prove that $\cv_2(\eta)$ diverges in the limit
$\eta\to 0$ when $\alpha<1/4$. Now, the integral $\cv_1(\eta)$ is the one which appears in the computations concerning the analytic fBm $\Gamma$, while the additional integral $\cv_2(\eta)$ is needed in order to handle the case of the real-valued fBm $\Re\Gamma$.
This fact had already been noted in \cite{Un}, where ${\cv}_1(\eta)$ (as part of the calculations
needed to compute ${\ca}_{ts}^{\eta}$) is evaluated in closed form involving Gauss' hypergeometric
function (see proof of Theorem 4.4 in \cite{Un}) -- in fact (see \cite{Un}, formula (4.36))
\begin{equation} {\cv}_1(\eta)\overset{\eta\to 0}{\to} \frac{\alpha(2\alpha-1)}{4\cos^2 \pi\alpha}
\ . \ \left[ \frac{2{\Gamma}(2\alpha-1){\Gamma}(2\alpha+1)}{{\Gamma}(4\alpha+1)}
+\frac{\cos 2\pi\alpha}{(2\alpha-1)(4\alpha-1)} \right] \ |t-s|^{4\alpha} \end{equation}
(a regular expression when $\alpha\to 1/4$ or $1/2$ as Taylor's formula proves). In the same article,
more general iterated integrals of the process ${\Gamma}^{\eta}$ are introduced {\it en passant}
 under the name of
{\em analytic iterated integrals} and shown to converge in the limit $\eta\to 0$; we reproduce these crucial results here. On the other hand, singularities of non-analytic iterated integrals
for $\alpha<1/4$ are analyzed in great details in \cite{Unt08}.

\smallskip

Here is how our article is structured. The first three sections are devoted to show how to solve differential equations in the complex plane like (\ref{eq:SDE}), under suitable assumptions on $b,\si$ and on the iterated integrals $X^k$: in Section \ref{sec:alg-integration}, we recall the basic features of the algebraic integration theory; Section \ref{sec:rough-2} aims at giving some details for the resolution of equation (\ref{eq:SDE}) in case of a rough signal with Hölder regularity $\ga>1/3$, while Section \ref{sec:rough-general} generalizes our considerations to an arbitrary Hölder regularity exponent in $(0,1/2)$. The remaining sections deal with the application of the general theory to the analytic fBm $\gga$: Section \ref{sec:def-fBm} is concerned with  the definition of this process,
 Section \ref{sec:grr}  with the proof of some general regularity results  for increments.  
Some useful (yet elementary) complex analysis preliminaries are given in Section
 \ref{sec:complex-analysis}. We then proceed to prove  the convergence of our approximations based on $\gga^\eta$: 
Section \ref{sec:cvgce-B-ep} deals with  $\gga^\eta$ itself, Section \ref{sec:levy-area} handles the case of 
the Levy area, while the general multiple integral case is treated in Section \ref{sec:multi-estim}.

\smallskip

\noindent
{\bf Notations:} Starting from Section \ref{sec:rough-2} and throughout the paper, the following notations concerning processes will be used. A generic $\ga$-Hölder function will be denoted by $X$. The analytic fBm defined on the complex upper half-plane is written $\Gamma=\{\Gamma_t;\, t\in\bar\Pi^+\}$, and its smooth approximation is denoted by $\Gamma^\ep$
or $\Gamma^\eta$. If $s,t\in\bar{\Pi}^+$, then $[s,t]=\{\lambda s+(1-\lambda)t\ |\ \lambda\in[0,1]\}\subset
\bar{\Pi}^+$ is the segment between $s$ and $t$. Generally speaking, $\Omega$ will denote a {\it bounded}
 neighborhood of $0$
in the closure of the upper half-plane $\bar{\Pi}^+$. Since this notation is generally used for probability spaces, we 
shall call $(\cu,\cf,\bp)$ the probability space under consideration here. Here is also a convention which will be used throughout the paper: for two real
 positive numbers, the relation $a\lesssim b$ stands for $a\le C \,b$, where $C$ is a given universal constant
(possibly depending continuously on $\alpha\in(0,1)$).

\section{Algebraic integration}
\label{sec:alg-integration}

Algebraic integration theory is conceived as an alternative to the popular rough paths
 analysis, and aims at solving differential equations driven by irregular processes with a
 minimal theoretical apparatus. Introduced in \cite{Gu} for a Hölder regularity of the driving noise
 $\ga>1/3$, it has then be extended to  arbitrary $\ga>0$ in  a  quite general setting (far beyond 
the geometric case) in  \cite{Gu2}. See also \cite{TT} for a detailed study of the case $\ga>1/4$. We have decided to recall some aspects of this formalism here for two main reasons: first, the integration theory we shall use takes place naturally in the upper complex plane $\Pi^+$. This induces some slight changes in the original setting, which we have chosen to outline. Second, we are able to deal with geometric rough paths in the current paper, which leads some simplifications in the analysis of the generalized integrals, compared for instance with \cite{Gu2}. We shall thus recall the main features of  algebraic integration in our context. This will also
 hopefully help to clarify the main assumptions which shall be checked on the process $\gga$.

\subsection{Increments}\label{sec:incr}

The extended integral we deal
with is based on the notion of increment, together with an
elementary operator $\der$ acting on them. These first notions are specifically introduced  in \cite{Gu,GT}, and we shall merely recall here their definition in the complex plane context. Consider an arbitrary neighborhood $\oom$ of 0 in the closure of the upper half-plane $\bar\Pi^+=\{z\in\C\ | \Im z\ge 0\}$. Then, for  a complex vector space $V$, and an integer $k\ge 1$, we denote by $\cac_k(\Omega;V)$ the set of functions $g : \oom^{k} \to V$ such
that $g_{t_1 \cdots t_{k}} = 0$
whenever $t_i = t_{i+1}$ for some $i\le k-1$.
Such a function will be called a
\emph{$(k-1)$-increment}, and we shall
set $\cac_*(V)=\cup_{k\ge 1}\cac_k(\Omega;V)$. The operator $\der$
alluded to above can be seen as an operator acting on
$k$-increments, 
and is defined as follows on $\cac_k(\Omega;V)$:
\begin{equation}
  \label{eq:coboundary}
\delta : \cac_k(\Omega;V) \to \cac_{k+1}(\Omega;V), \qquad
(\delta g)_{t_1 \cdots t_{k+1}} = \sum_{i=1}^{k+1} (-1)^i g_{t_1
  \cdots \hat t_i \cdots t_{k+1}} ,
\end{equation}
where $\hat t_i$ means that this particular argument is omitted.
Then a fundamental property of $\der$, which is easily verified,
is that
$\delta \delta = 0$, where $\delta \delta$ is considered as an operator
from $\cac_k(\Omega;V)$ to $\cac_{k+2}(\Omega;V)$, so $(\cac_*(\Omega;V),\delta)$ is a cochain complex.
 we shall denote $\cZ\cac_k(\Omega;V) = \cac_k(\Omega;V) \cap \text{Ker}\delta$
and $\cb \cac_k(\Omega;V) =\cac_k(\Omega;V) \cap \text{Im}\delta$.

\smallskip

Some simple examples of actions of $\der$,
which will be the ones we shall really use throughout the paper,
 are obtained by letting
$g\in\cac_1$ and $h\in\cac_2$. Then, for any $t,u,s\in\oom$, we have
\begin{equation}
\label{eq:simple_application}
  (\der g)_{ts} = g_t - g_s,
\quad\mbox{ and }\quad
(\der h)_{tus} = h_{ts}-h_{tu}-h_{us}.
\end{equation}
Furthermore, it is readily checked that
the complex $(\cac_*,\delta)$ is \emph{acyclic}, i.e.
$\cZ \cac_{k}(\Omega;V) = \cb \cac_{k}(\Omega;V)$ for any $k\ge 1$. 

\smallskip

Let us mention at this point some conventions on products of increments which will be used in the sequel: assuming for the moment
that $V=\C$, set $\cac_k(\Omega;\C)=\cac_k(\Omega)$. Then
the complex $(\cac_*(\Omega),\delta)$ is an (associative, non-commutative)
graded algebra once endowed with the following product:
for  $g\in\cac_n(\Omega) $ and $h\in\cac_m(\Omega) $ let  $gh \in \cac_{n+m-1}(\Omega) $
be the element defined by
\begin{equation}\label{cvpdt}
(gh)_{t_1,\dots,t_{m+n-1}}=g_{t_1,\dots,t_{n}} h_{t_{n},\dots,t_{m+n-1}},
\quad
t_1,\dots,t_{m+n-1}\in\oom.
\end{equation}
The pointwise multiplication of $g,\hat g\in\cac_n(\Omega)$, denoted by $g\circ\hat g$, is also defined by:
$$
(g\circ\hat g)_{t_1,\dots,t_{n}}= g_{t_1,\dots,t_{n}} \, \hat g_{t_1,\dots,t_{n}},
\quad
t_1,\dots,t_{n}\in\oom.
$$

\smallskip

Our future discussions will mainly rely on
$k$-increments with $k \le 2$, for which we shall use some
analytic assumptions. Namely, sticking to the case $V=\C^d$ for $d\ge 1$,
we measure the size of these increments by H\"older norms
defined in the following way: for $f \in \cac_2(\Omega;V)$ let
$$
\norm{f}_{\mu} \equiv
\sup_{s,t\in\oom}\frac{|f_{ts}|}{|t-s|^\mu},
\quad\mbox{and}\quad
\cac_2^\mu(V)=\lcl f \in \cac_2(\Omega;V);\, \norm{f}_{\mu}<\infty  \rcl.
$$
In the same way, for $h \in \cac_3(\Omega;V)$, set 
\begin{eqnarray}
  \label{eq:normOCC2}
  \norm{h}_{\ga,\rho} &=& \sup_{s,u,t\in\oom}
\frac{|h_{tus}|}{|u-s|^\ga |t-u|^\rho}\\
\norm{h}_\mu &\equiv &
\inf\left \{\sum_i \norm{h_i}_{\rho_i,\mu-\rho_i} ;\, h  =\sum_i h_i,\, 0 < \rho_i < \mu \right\} ,\nonumber
\end{eqnarray}
where the last infimum is taken over all sequences $\{h_i \in \cac_3(\Omega;V) \}$ such that $h
= \sum_i h_i$ and for all choices of the numbers $\rho_i \in (0,z)$.
Then  $\norm{\cdot}_\mu$ is easily seen to be a norm on $\cac_3(\Omega;V)$, and we set
$$
\cac_3^\mu(V):=\lcl h\in\cac_3(\Omega;V);\, \norm{h}_\mu<\infty \rcl.
$$
Notice that, in order to avoid ambiguities, we shall use the notation $\cn[\cdot\,;\cac_k^\mu(\Omega;V)]$, 
instead of $\norm{\cdot}_\mu$, to denote the Hölder norms in the spaces $\cac_k(\Omega;V)$.
Finally,
let $\cac_3^{1+}(\Omega;V) = \cup_{\mu > 1} \cac_3^\mu(\Omega;V)$,
and remark that the same kind of norms can be considered on the
spaces $\cZ \cac_3(\Omega;V)$, leading to the definition of normed subspaces
$\cZ \cac_3^\mu(\Omega;V)=\cac_3^{\mu}(\Omega;V)\cap \cZ \cac_3(\Omega;$ $V)$ and
 $\cZ \cac_3^{1+}(\Omega;V)=\cup_{\mu>1} \cac_3^{\mu}(\Omega;V)\cap \cZ \cac_3(\Omega;V)$.
 It also turns out to be useful to consider the  spaces of continuous
$(k-1)$ increments $\cac_k^0(\Omega;V)$, equipped with the norm
 $\cn[h;\, \cac_k^0(\oom;V)]=\sup_{t_1,\ldots,t_k\in\oom}|h_{t_1,\ldots,t_k}|$. Let us mention once and for all that all the Hölder spaces  we are considering in this article are complete.

\subsection{Iterated integrals on the upper-half plane}\label{cpss}

The iterated integrals of analytic functions on $\Omega\cap\Pi^{+}$ are 
particular cases of elements of $\cac_*$ which will be of interest for
us. Let us recall  some basic  rules for these objects. Set $\cac_k^{{\om}}(\Omega\cap\Pi^{+})
=\cac_k^{{\om}}(\Omega\cap\Pi^{+};\,\C)$ for the set of analytic $(k-1)$-increments from $\Omega\cap\Pi^{+}$ to $\C$, and 
consider $f,g\in\cac_1^{{\om}} (\Omega\cap\Pi^{+})$. Then the integral $\int dg \,
f$, which will also be denoted by
$\cj(dg \,  f)$, can be considered as an element of
$\cac_2^{{\om}}(\Omega\cap\Pi^{+})\equiv\cac_2^{{\om}}(\Omega\cap\Pi^{+};\C)$. That is, for $s,t\in\Omega\cap\Pi^{+}$, we set
$$
\cj_{ts}(dg \  f)=
\left(\int  dg f \right)_{ts} = \int_{[s,t]}  dg_u f_u.
$$
The multiple integrals can also be defined in the following way:
given a smooth element $h \in \cac_2^{{\om}}$ and $s,t\in\Omega\cap\Pi^{+}$, we set
$$
\cj_{ts}(dg\ h )\equiv
\left(\int dg h \right)_{ts} = \int_{[s,t]} dg_u h_{us} .
$$
In particular, the double integral $\cj_{ts}( df^3df^2\,f^1)$ is defined, for
$f^1,f^2,f^3\in\cac_1^{{\om}}(\Omega\cap\Pi^{+})$, as
$$
\cj_{ts}( df^3df^2\,f^1)
=\lp \int df^3df^2\,f^1  \rp_{ts}
= \int_{[s,t]} df_u^3 \,\cj_{us}\lp  df^2 \, f^1 \rp .
$$
Now, suppose that the $n\textsuperscript{th}$ order iterated integral of $df^n\cdots df^2 \,f^1$, still denoted by $\cj(df^n$ $\cdots df^2 \,f^1)$, has been defined for
$f^1,f^2\ldots, f^n\in\cac_1^{{\om}}(\Omega\cap\Pi^{+})$.
Then, if $f^{n+1}\in\cac_1^{{\om}}(\Omega\cap\Pi^{+})$, we set
\begin{equation}\label{multintg}
\cj_{ts}(df^{n+1}df^n \cdots df^2 f^1)=
\int_{[s,t]}  df_u^{n+1}\, \cj_{us}\lp df^n\cdots df^2 \,f^1\rp\,,
\end{equation}
which defines the iterated integrals of smooth functions recursively.
Observe that a $n\textsuperscript{th}$ order integral $\cj(df^n\cdots df^2 df^1)$, where we have simply replaced $f^1$ by $df^1$, could be defined along the same lines.

\smallskip

The following relations between multiple integrals and the operator $\der$ are easily checked  for analytic functions. They are also a prototype of the algebraic relations we shall impose in the rough setting:
\begin{proposition}\label{prop:dissec}
Let $f,g$ be two elements of $\cac_1^{{\om}}(\Omega\cap\Pi^{+})$. Then, recalling the convention
(\ref{cvpdt}), it holds that
$$
\der f = \cj( df), \qquad
\der\lp \cj( dg f)\rp = 0, \qquad
\der\lp \cj (dg df)\rp =  (\der g) (\der f) = \cj(dg) \cj(df),
$$
and, in general,
$$
 \der \lp \cj( df^n \cdots df^1)\rp   =  \sum_{i=1}^{n-1}
\cj\lp df^n \cdots df^{i+1}\rp \cj\lp df^{i}\cdots df^1\rp
$$
which is simply a way of rewriting the multiplicative property (\ref{eq:*}).
\end{proposition}


\subsection{The complex sewing map}

The algebraic integration theory heavily relies on a generalization of Young integrals to increments in $\cac_2$ (as explained in \cite{Gu,GT}). This generalization is obtained through a map which is often called the \emph{sewing map} according to the terminology of \cite{FP}. We present here a construction of this map adapted to the complex plane context, for which some additional notation is needed: for $j\ge 1$ and $\beta>0$, let $\cac_j^{{\rm m},\beta}(\Omega;V)$,
 where ${\rm m}$ stands for \emph{multiparametric}, be the subspace of $\cac_j(\Omega;V)$ induced by the semi-norm:
\begin{equation}\label{eq:13}
\cn[h;\,\cac_j^{{\rm m},\beta}(\Omega;V)]=
\sup\lcl
\frac{| h_{s_1+\imath\ep,\ldots, s_j+\imath\ep}-h_{s_1,\ldots ,s_j}|}{\ep^{\beta}};\,
\ep\in[0,1],\ s_1,\ldots,s_j\in\oom
\rcl.
\end{equation}
The definition is better understood when one thinks of the  regularity properties of the boundary value of hyperfunctions (see
\cite{Mor93}, Theorem 3.9.8); for instance it is known that if $f$ is analytic on $\Pi^+$, and
 $\sup_{x\in\R}|f(x+\imath y)|\lesssim y^{-\beta}$ for 
some $\beta$, then the boundary value of $f$ is a distribution.

In order to define a sewing map on the complex plane, we shall then use the following norm defined for  $(j-1)$-increments:
$$
\cn_{\C}[h;\,\cac_j^{\mu,\beta}(\Omega;V)]= \cn[h;\,\cac_j^{\mu}(\Omega;V)]
+\cn[h;\,\cac_j^{{\rm m},\beta}(\Omega;V)],
$$
and as usual now, we call $\cac_j^{\mu,\beta}(\Omega;V)$ the  associated normed space.

\begin{proposition}[{\bf construction of the complex sewing map $\Lambda$}]  \label{prop:Lambda}
Let $\mu >1$, $\beta>0$, and assume that $h\in\cac_3(\Omega;V)$ satisfies the following hypotheses:

\smallskip

\noindent
(i) $h$ is an element of $\cz \cac_3^\mu(\Omega;V)\cap\cac_3^{{\rm m},\beta}(\Omega;V)$. In particular, $\cn_{\C}[h;\,\cac_3^{\mu,\beta}(\Omega;V)]$ is finite.

\smallskip

\noindent
(ii)  For $t,u,s\in\Pi^+$, $h_{tus}$ can be written as $h_{tus}=[\der(\cj(df \, r))]_{tus}$ for an analytic function $f$ and an increment $r\in\cac_2$ such that, for any $s\in\Pi^+$, the function $u\mapsto r_{us}$ is analytic.

\smallskip

\noindent
Then, for any $1<\nu<\mu$, there exists a unique $\Lambda h \in \cac_2^\nu(\Omega;V)$ such that $\der( \Lambda h )=h$. Furthermore, there exists a strictly positive constant $c_\nu$ such that
\begin{eqnarray} \label{contraction}
\cn[\Lambda(h);\, \cac_2^\nu(\Omega;V)] \leq c_\nu \, \cn_{\C}[h;\,\cac_3^{\mu,\beta}(\Omega;V)].
\end{eqnarray}
Calling $\ca_3^{\mu,\beta}(\Omega;V)$ the set of increments satisfying conditions (i) and (ii) above, this gives rise to a  continuous linear
map $\laa:  \ca_3^{\mu,\beta}( \Omega;V) \rightarrow \cac_2^\nu(\Omega;V)$ such that $\der \laa =\id_{ \ca_3^{\mu,\beta}( \Omega;V)}$.
\end{proposition}

\begin{proof}
The proof of this proposition is an extension of \cite[Proposition 2.3]{GT}, for which we refer for further details. In particular, the uniqueness of $g:=\laa(h)$ can be proven just as in the above reference. We shall thus outline the proof of the existence part of the increment $g$.

\smallskip

The increment  $g$ can be defined in a natural way in two cases: (1) When $t,s\in\oom\cap\R$, then $g_{ts}$ can be constructed as in the proof of \cite[Proposition 2.3]{GT}. (2) When $t,s\in\oom\cap\Pi^+$, then $g_{ts}$ can simply be defined as $\cj_{ts}(df \, r)$. By construction, we thus have that $[\der g]_{tus}=h_{tus}$ when $t,u,s\in\oom\cap\R$ or $t,u,s\in\oom\cap\Pi^+$.

\smallskip

In order to define $g_{ts}$ for $t\in\Pi^+$ and $s\in\R$, we use two limiting procedures: first, for $\ep>0$, we set $g^{1,\ep}_{ts}=\cj_{t+\imath\ep,s+\imath\ep}(df \, r)$, which is defined in the Riemann sense. Next, another candidate for an approximating sequence of $g$ is $g^{2,\ep}$ constructed along the same lines as in \cite[Proposition 2.3]{GT}: we know
that there exists a $B\in\cac_2$ such that $\der B=h$.
For $n\ge 0$, consider the
dyadic partition $\{r^{\ep,n}_j; 0\le j\le 2^n  \}$ of the segment $[s+\imath\ep,t]$, where
\begin{equation*}
r^{\ep,n}_j=s+\imath\ep+\frac{(t-s)j}{2^n},
\quad\mbox{for}\quad 0\le j\le 2^n.
\end{equation*}
Then, for $n\ge 0$  set $M_{ts}^{\ep,n}=B_{ts}-\sum_{j=0}^{2^n-1} B_{r^{\ep,n}_{j+1},r^{\ep,n}_{j}}$.
In this context, it can be shown that 
\begin{equation}\label{eq:19}
M_{ts}^{\ep,n+1}-M_{ts}^{\ep,n}=\sum_{j=0}^{2^{n}-1} h_{r^{\ep,n+1}_{2j+2},r^{\ep,n+1}_{2j+1},r^{\ep,n+1}_{2j}},
 \end{equation}
from which the convergence of $M_{ts}^{\ep,n}$ is easily deduced. We call $g^{2,\ep}:=\lim_{n\to\infty}M^{\ep,n}$. By a uniqueness argument, restricted on the segment $[s+\imath\ep,t+\imath\ep]$ and elaborated as in the real case of \cite[Proposition 2.3]{GT}, it can also be shown that $g^{1,\ep}=g^{2,\ep}$. 

\smallskip

Furthermore, we claim that for any $1<\nu<\mu$ and $\ep,\eta$ small enough, we have 
\beq\label{eq:28}
\cn[g^{2,\ep}-g^{2,\eta};\, \cac_2^{\nu}(\Omega;V)] \lesssim
 \cn_{\C}[h;\,\cac_3^{\mu,\beta}(\Omega;V)] |\ep-\eta|^{\hat\beta},
\eeq
for a certain $\hat\beta>0$. Indeed, we have defined $g^{2,\ep}$ as $\sum_{n\ge 0} M^{\ep,n+1}-M^{\ep,n}$, and it is thus sufficient to show that
\beq\label{eq:29}
|\Delta_{ts}^{n}(\ep,\eta)|
\lesssim \cn_{\C}[h;\,\cac_3^{\mu,\beta}(\Omega;V)] |\ep-\eta|^{\hat\beta} \, 2^{-n \hat\ga} \,
|t-s|^{\nu}
\eeq
for some $\hat\ga>0$, where we have set $\Delta_{ts}^{n}(\ep,\eta):=(M_{ts}^{\ep,n+1}-M_{ts}^{\ep,n})- (M_{ts}^{\eta,n+1}-M_{ts}^{\eta,n})$. Now, thanks to our decomposition (\ref{eq:19}), it is readily checked for an arbitrary $\rho\in(0,1)$ that
\begin{multline*}
|\Delta_{ts}^{n}(\ep,\eta)| \le
\sum_{j=0}^{2^n-1} 
\lln h_{r^{\ep,n+1}_{2j+2},r^{\ep,n+1}_{2j+1},r^{\ep,n+1}_{2j}}
-h_{r^{\eta,n+1}_{2j+2},r^{\eta,n+1}_{2j+1},r^{\eta,n+1}_{2j}}\rrn^{1-\rho}  \\
\times \lp\lln h_{r^{\ep,n+1}_{2j+2},r^{\ep,n+1}_{2j+1},r^{\ep,n+1}_{2j}} \rrn
+\lln h_{r^{\eta,n+1}_{2j+2},r^{\eta,n+1}_{2j+1},r^{\eta,n+1}_{2j}}\rrn\rp^{\rho},
\end{multline*}
and hence, invoking the fact that $h$ is an element of $\cac_3^{{\rm m},\beta}(\Omega;V)$ and $\cac_3^\mu(\Omega;V)$, we obtain:
$$
|\Delta_{ts}^{n}(\ep,\eta)| \le 3 \,
\cn_{\C}[h;\,\cac_3^{\mu,\beta}(\Omega;V)] \, |\ep-\eta|^{(1-\rho)\beta} \, 2^{-n (\rho\mu-1)} \,
|t-s|^{\rho\mu},
$$
which immediately entails (\ref{eq:29}) for $\rho\mu=\nu$, and thus (\ref{eq:28}), since $\rho$ can be taken arbitrarily close to 1. This  implies the convergence of $g^{2,\ep}$ towards an increment $g$ as $\ep\to 0$, and allows to define $g_{ts}$ for $t\in\Pi^+$ and $s\in\R$. Notice that $g=\lim_{\ep\to 0}g^{1,\ep}=\lim_{\ep\to 0}g^{2,\ep}$, the limit being understood in $\cac_2^\nu(\oom;V)$.

\smallskip

We now have to check that $(\der g)_{tus}=h_{tus}$ in two remaining cases: (1) When $t,u\in\Pi^+$ and $s\in\R$.  (2) When $t\in\Pi^+$ and $u,s\in\R$. In fact, the two cases can be treated similarly, and we focus on the first one. To this purpose, we resort to the approximation $g_{ts}=\lim_{\ep\to 0}g^{1,\ep}_{ts}$. This yields the relation:
$$
\der g_{tus}=\lim_{\ep\to 0} \der\lc  \cj(df\, r)\rc_{t+\imath\ep,u+\imath\ep,s+\imath\ep}
=\lim_{\ep\to 0}  h_{t+\imath\ep,u+\imath\ep,s+\imath\ep}=h_{tus},
$$
where the second equality is obtained thanks to our assumption (ii) on the increment $h$. We have thus proved our claim $\der g=h$. It remains to prove that $g\in\cac_3^\nu$, which can be done exactly as in \cite[Proposition 2.3]{GT}.

\end{proof}

It is also worth mentioning that the operator $\Lambda$ (called {\em complex sewing map}) can be related to some Riemann type sums, which is a way to link the objects constructed so far with a generalized notion of integral:
\begin{proposition}[Integration of small increments]
\label{prop:integration}
For any 1-increment $g\in\cac_2 (\oom;V)$ such that $\der g$ satisfies the assumptions of Proposition \ref{prop:Lambda}, set $\delta f = (\id-\Lambda \delta) g$. Then
$$
(\delta f)_{ts} = \lim_{|D_{ts}| \to 0} \sum_{i=0}^n g_{t_{i+1}\, t_i},
$$
where the limit is over any partition $D_{ts} = \{t_0=t,\dots,
t_n=s\}$ of the segment $[s,t]\subset\bar\Pi^+$, whose mesh tends to zero. The
1-increment $\delta f$ is thus the indefinite integral of the 1-increment $g$.
\end{proposition}


\section{Rough path analysis of order 2}\label{sec:rough-2}
As we mentioned before, our complex setting leads us naturally to solve differential equations in a neighborhood of 0 in the upper half-plane $\bar\Pi^+$. This framework being rather unusual, we shall try to explain its main features on the simplest non-trivial example one can think of, namely the case of a driving signal with Hölder continuity exponent greater than $1/3$. More specifically, we assume in this section that $X$ satisfies the following condition:
\begin{hypothesis}\label{hyp:1}
The signal $X$ is an element of $\cac^\ga_1(\oom;\C^d)$ with $1/3<\ga\le 1/2$.  Furthermore, its restriction to $\Pi^+$ is analytic.
\end{hypothesis}
Under these assumptions, we shall now explain how to solve equation (\ref{eq:SDE}). Once again, these considerations follow closely the methodology introduced in \cite{Gu}, but we include them here because of the slight changes due to our complex plane situation.

\smallskip

In the sequel of the paper, we shall use indistinctly $\int_{[s,t]}  dg \, f$
or $\cj_{ts}(dg \, f)$ for the integral of a function $f$ with
respect to a given increment $dg$ on the segment $[s,t]\subset\oom$.
Observe that the second notation $\cj_{ts}(dg \, f)$ aims at avoiding
some cumbersome notations in our future computations. Recall also
that we wish to solve an equation of the form 
\begin{equation}
\label{red:diff-eq}
y_t=a+ \int_{[0,t]}  dX^*_s \, \si(y_s),
\quad t\in\oom,
\end{equation}
where $X$ is a given $\ga$-H\"older continuous path from $\oom$ to $\R^{d}$,
with $1/3<\ga<1/2$, and $X^*$ is the transpose of the path  $X$,
 considered as a $\C^{1,d}$-valued process. Notice that in the last equation, we have chosen to use the slightly unusual convention of multiplying the coefficient $\si$ by the driving process $X$ in order to simplify a little our further expansions. With respect to equation (\ref{eq:SDE}), we have also chosen to skip the drift term $b$ for notational sake, though the inclusion of such a drift term would be technically easy.

\vspace{0.3cm}

Before going into the technical resolution of (\ref{red:diff-eq}),
let us make some heuristic considerations about the form that
a candidate solution should have: set 
$\hsi_t=\si(y_t)$, and suppose we have
been able to exhibit a solution $y$ to (\ref{red:diff-eq}),
such that $y\in\cac_1^\ka$ for a given $1/3<\ka<\ga$.
Then the integral form of our equation, for $t\in\oom$, can be 
read as
\begin{equation}\label{intg:delay}
y_t=a+\int_{[0,t]}  dX^*_u \, \hsi_u.
\end{equation}
Our approach to generalized integrals induces us to work, instead of 
(\ref{intg:delay}),  with increments
of the form $(\delta y)_{ts}=y_t-y_s$. However, it is easily checked that,
provided one is  given a reasonable notion of integral, one can decompose
(\ref{intg:delay}) into
$$
(\delta y)_{ts}=\int_{[s,t]} dX^*_u \, \hsi_u = (\delta X^*)_{ts} \, \hsi_s+\rho_{ts},
\quad\mbox{ with }\quad
\rho_{ts}=\ist dX^*_u \, (\hsi_u-\hsi_s) .
$$
We have thus obtained a decomposition of $y$ of the form
$\delta y=\delta X^*\,\hsi+\rho$. Let us see, still at a heuristic level,
what is the regularity one can expect on $\hsi$ and $r$: first, suppose that
$\si$ is an analytic function, and  assume that $y_t\in B(0,M)$ for all $t\in\oom$,
 where $B(0,M)=\{z\in\C;\, |z|\le M\}$. Then $\hsi$ is bounded and
$$
|\hsi_t-\hsi_s|
\le
\norm{\si'}_{\infty,M}
\lln  y_{t}-y_{s}  \rrn
\le
\norm{\si'}_{\infty,M} \,\cn[y;\, \cac_1^\ga(\Omega)] \,|t-s|^{\ga},
$$
where the notation $\cn$ has been introduced in Section \ref{sec:incr} and where $\|\si'\|_{\infty,M}$ stands
 for $\sup_{z\in B(0,M)}|\si'(z)|$. Thus, still with the notations of Section \ref{sec:incr},
we have that $\hsi\in\cac_1^\ga(\Omega)\cap\cac_1^0(\Omega)$. As far as $\rho$ is
concerned, provided one can define the integral
$\ist  dX^*_u\, (\hsi_u-\hsi_s)=\ist dX^*_u\, (\delta\hsi)_{su}$, it should
inherit both regularities of $\delta\hsi$ and $dx$. Thus, one should
expect that $\rho\in\cac_2^{2\ka}$. In conclusion, we have found that
$\delta y$ should be decomposable into
\begin{equation}\label{first:structure}
\delta y = \delta X^*\, \hsi + \rho,
\quad\mbox{ with }\quad
\hsi\in\cac_1^\ga(\Omega)\cap\cac_1^0(\Omega)
\mbox{ and }
\rho\in\cac_2^{2\ka}(\Omega).
\end{equation}
This motivates the definition of what we call the class of {\it analytic controlled paths} of order 1,
in which we shall solve equation (\ref{red:diff-eq}) when $\ga>1/3$:
\begin{definition}[{\bf space of analytic controlled paths $\cq_{\ka,a}(\Omega)$}]
\label{def39}
Let $z$ be a process in $\cac_1^\ka(\Omega;\C^{k})$ with $1/3<\ka\le\ga$.
We say that $z$ is an analytic controlled path of order 1 based on $X$, if
$z_0=a$, which is a given initial condition in $\C^{k}$,
and for $i=1,\ldots,k$, the increment $\der z(i)\in\cac_2^\ka(\Omega;\C)$ can be decomposed into
\begin{equation}\label{weak:dcp}
\der z(i)=\der X(j)\, \zeta(j,i) + r(i),
\quad\mbox{i. e.}\quad
(\der z)_{ts}(i)=(\der X)_{ts}(j)\, \zeta_s(j,i)  + r_{ts}(i),
\quad s,t\in\oom,
\end{equation}
with $\zeta\in\cac_1^\ka(\Omega;\C^{d,k})$ and $r$ is a regular part
such that $r\in\cac_2^{2\ka}(\Omega;\C^{k})$. In our complex plane setting, we assume moreover that $z,\zeta$ and $u\mapsto r_{us}$ are analytic paths when restricted to $\Pi^+$, for any $s\in\Pi^+$.
The space of analytic controlled
paths will be denoted by $\cq_{\ka,a}(\Omega;\C^{k})$, and a process
$z\in\cq_{\ka,a}(\Omega;\C^{k})$ can be considered in fact as a couple
$(z,\zeta)$. The natural semi-norm on $\cq_{\ka,a}(\Omega;\C^{k})$ is given
by
\begin{multline*}
\cn[z;\cq_{\ka,a}(\Omega;\C^{k})]=
\cn[z;\cac_1^{\ka}(\Omega;\C^{k})]
+ \cn[\zeta;\cac_1^{0}(\Omega;\C^{d,k})]
+ \cn[\zeta;\cac_1^{\ka}(\Omega;\C^{d,k})]\\
+\cn[r;\cac_2^{2\ka}(\Omega;\C^{k})],
\end{multline*}
with 
$\cn[\zeta;\cac_1^{0}(V)]=\sup_{0\le s\le T}|\zeta_s|_V$. 
\end{definition}

\vspace{0.3cm}

With this definition at  hand, here is the global strategy 
we shall adopt in order to solve equation (\ref{red:diff-eq}):
\begin{enumerate}
\item
Study the stability of $\cq_{\ka,a}(\Omega;\C^{k})$ under an analytic map
$\vp:\C^{k}\to\C^{n}$.
\item
Define rigorously the integral $\int dX^*_u\,z_u =\cj(dX^*\,z )$
for an analytic controlled path $z$ and compute its decomposition
(\ref{weak:dcp}).
\item
Solve equation (\ref{red:diff-eq}) in the space $\cq_{\ka,a}(\Omega;\C^{k})$
by a fixed point argument.
\end{enumerate}
Let us go on with this program, and first see how smooth functions act on analytic controlled paths.


\subsection{Action of analytic maps on controlled paths}


The action of a smooth function on an analytic controlled path can be 
summarized in the following proposition, for which we need an additional notation: given an analytic function $\vp :\C^{k}\to\C^{n}$, we set $\vp^{i}(z)$ for the $i\textsuperscript{th}$ coordinate of $\vp(z)$ and $\partial_j\vp^{i}$ for the derivative of $\vp^{i}$ with respect to the $j\textsuperscript{th}$ coordinate of $z$.
\begin{proposition}\label{cp:weak-phi}
Let $z\in\cq_{\ka,a}(\Omega;\C^{k})$ with decomposition (\ref{weak:dcp}). Let 
$\vp :\C^{k}\to\C^{n}$ be  an analytic function, and set $\hz=\vp(z)$, $\ha=\vp(a)$.
Then $\hz\in\cq_{\ka,\ha}(\Omega;\C^{n})$, and it can be decomposed into
$$
\der \hz(i)= \der X(k) \hat\zeta(k,i)  +\hr(i), \quad
i=1,\ldots,n
$$
with
$$
\hat\zeta(k,i)= \zeta(k,j)\partial_j\vp^{i}(z)
\quad\mbox{ and }\quad
\hr(i)= r(j) \partial_j\vp^{i}(z) + \lc \der(\vp^{i}(z))-\der z(j)\partial_j\vp^{i}(z) \rc.
$$
Furthermore, assuming that $z_s\in B(0,M)$ if $s\in\oom$, where $B(0,M)$ stands for the set $\{z\in\C^n;\, |z|\le M\}$, the following bound holds true:
\begin{equation}\label{bnd:phi}
\cn[\hz;\cq_{\ka,\ha}(\Omega;\C^{n})]\le
c_{M}\lp \cn[z;\cq_{\ka,a}(\Omega;\C^{n})]+\cn^2[z;\cq_{\ka,a}(\Omega;\C^{n})]  \rp.
\end{equation}
\end{proposition}
\begin{proof}
The algebraic part of the assertion is quite straightforward. Just write
$$
(\der\hz)_{ts}(i)=\vp^{i}(z_t)-\vp^{i}(z_s)=
(\der z)_{ts}(j) \partial_j\vp^{i}(z_s)
+ \vp^i(z_t)-\vp^i(z_s)-(\der z)_{ts}(j) \partial_j\vp^{i}(z_s),
$$
and plugging expression (\ref{weak:dcp}) instead of $\der z(j)$ above, our first assertion is easily shown.

\smallskip

In order to give an estimate for $\cn[\hz;\cq_{\ka,\ha}(\Omega;\C^{n})]$, one must of  course
 establish bounds for $\cn[\hz;\cac_1^{\ka}(\Omega;\C^{n})]$,
 $\cn[\hat\zeta;\cac_1^{0}(\Omega;\C^{n,k})]$,  $\cn[\hat\zeta;\cac_1^{\ka}(\Omega;\C^{n,k})]$
and $\cn[\hr;\cac_2^{2\ka}(\Omega;\C^{n})]$. Let us focus on the last of these
estimates, the other ones being quite similar. First notice that
$\hr=\hr^1+\hr^2$ with
\begin{equation}\label{def:hr12}
\hr^1(i)=r(j) \partial_j\vp^{i}(z) 
\quad\mbox{ and }\quad
\hr^2(i)=\der(\vp^{i}(z))-\der z(j)\partial_j\vp^{i}(z).
\end{equation}
Now, since $\nabla\vp$ is a bounded $\C^{n,k}$-valued function on $B(0,M)$, we have
\begin{equation}\label{est:hr1}
\cn[\hr^1;\cac_2^{2\ka}(\Omega;\C^{n})]
\le
c_{1,M} \cn[r;\cac_2^{2\ka}(\Omega;\C^{k})].
\end{equation}
Moreover, with the same kind of argument for $\nabla^2\vp$, we also get:
\begin{equation*}
|\hr^2_{ts}|\le c_{2,M} |(\der z)_{ts}|^2
\le c_{2,M} \cn^{2}[z;\cac_1^\ka(\Omega;\C^{k})]|t-s|^{2\ka},
\end{equation*}
which yields
\begin{equation}\label{est:hr2}
\cn[\hr^2;\cac_2^{2\ka}(\Omega;\C^{n})]\le
c_{2,M} \cn^2[r;\cac_2^{2\ka}(\Omega;\C^{k})],
\end{equation}
and thus
we obtain
$$
\cn[\hr;\cac_2^{2\ka}(\Omega;\C^{n})]\le
c_{3,M} \lp 1+\cn^2[r;\cac_2^{2\ka}(\Omega;\C^{k})]\rp.
$$
The analytic assumptions on the increments $\hat z, \hat\zeta$ and $\hat r$ are then readily checked, and this concludes the proof.

\end{proof}

\subsection{Integration of analytic controlled paths}

The aim of this section is to define the integral $\cj(dX^*\,m)$
for an analytic controlled process $m\in\cq_{\ka,a}(\Omega;\C^{d})$, 
admitting the decomposition given by (\ref{weak:dcp}). Namely, we assume that $m$ can be decomposed as:
\begin{equation}\label{dcp:delta-m}
\der m(i)=\der X(j)\, \mu(j,i) + r(i), \quad i=1,\ldots,d,
\end{equation}
where $\mu \in\cac_1^\ka(\Omega;\C^{d,d})$ and $r\in\cac_2^{2\ka}(\Omega;\C^{d})$.
In order to see how the integral $\cj(dX^*\,m)$ may look like, let us treat first
the case of smooth processes $X,\mu$ and $r$, and see how 
$\cj(dX^*\,m)$ can be expressed in terms of the operators $\delta$
and $\laa$: in the regular case, $\cj(dX^*\,m)$ is well-defined, and
we have
$$
\ist  dX_u^* \, m_u = [X_t-X_s]^* \, m_s+ \ist  dX_u^* \, [m_u-m_s] 
$$  for $s, t\in\oom$,
or in other words
$$
\cj(dX^* \, m)= \der X^*\, m+ \cj(dX^*\, \der m).
$$
Let us plug the  decomposition (\ref{dcp:delta-m}) into this expression,
which yields
\begin{eqnarray}\label{eq:imdx}
\cj(dX^*\, m)&=& \der X^*\, m + \cj\lp dX(i) \, \der X(j)\, \mu(j,i) \rp + \cj(dX^*\, r).
\end{eqnarray}
Let us transform now the term $\cj(dX(i) \der X(j)\, \mu(j,i) )$: it is easily seen that
$$
\cj_{ts}\lp dX(i) \, \der X(j)\, \mu(j,i) \rp
= \ist   dX_u(i) \, \der X_{us}(j) \ \mu_s(j,i)
=\bx_{ts}^{\bf 2}(i,j) \, \mu_s(j,i) 
$$ 
for $s,t\in\oom$, where $\xd$ is the $\cac_2(\Omega;\C^{d,d})$-function defined by
$$
\xdst=\ist dX_u \otimes [\der X]_{su},
\quad\mbox{i.e.}\quad
[\xdst](i,j)= \ist dX_u(i) \, [\der X]_{su}(j),
\quad 1\le i,j\le d.
$$
Inserting this expression  into (\ref{eq:imdx}), we get
\begin{equation}\label{exp2:imdx}
\cj(dX^*\, m)= \der X^*\, m + \xd(i,j) \, \mu(j,i) + \cj(dX^*\, r).
\end{equation}

\smallskip

Let us focus now on the term $\xdst$: when $X$ is a smooth
process, it is readily checked that we have
$$
[\der\xd]_{tus}
=\bx_{ts}^{\bf 2}-\bx_{su}^{\bf 2}-\bx_{ut}^{\bf 2}
=(\der X)_{tu}\otimes(\der X)_{us}.
$$
This decomposition of $\der\xd$ into a product of increments is
the fundamental algebraic property we shall use in order to extend
the above integral to non-smooth processes, and thus, in the sequel 
of this section, we shall assume the following:
\begin{hypothesis}\label{hyp:x1}
The path $X$ is $\C^d$-valued
$\ga$-H\"older with $\ga>1/3$, it satisfies Hypothesis \ref{hyp:1}, and  it admits a so-called {\em L\'evy area},
that is, a process 
$\xd\in\cac_2^{2\ga}(\Omega;\C^{d,d})$, defined formally as $\xd=\cj(dX \otimes dX)$, and satisfying
$$
\der\xd=\der X\otimes \der X,
\quad\mbox{i.\!\! e.}\quad
\lc (\der\xd)_{tus} \rc(i,j)
=
[\der X]_{tu}(i) [\der X]_{us}(j),
$$
for any $t,u,s\in\oom$ and $i,j\in\{1,\ldots,d  \}$. We also suppose that $\xd$ belongs to the space $\cac_2^{{\rm m},\ga}(\oom;\C^{d,d})$ defined by relation (\ref{eq:13}).
Notice that under our  assumptions, the increment $\xd$ is analytic on $\Pi^+\times\Pi^+$.
\end{hypothesis}
Let us keep this hypothesis in mind, while we finish the analysis
of the smooth case: it remains to find a suitable expression for
$\cj(dX^*\, r)$. To this purpose, let us write (\ref{exp2:imdx})
as
\begin{equation}\label{exp1:irhodx}
\cj(dX^*\, r)=\cj(dX^*\, m)-\der X^*\, m - \xd \cdot \mu^*,
\end{equation}
where in the above expressions, we have denoted by $M\cdot N$ the Hilbert-Schmidt 
inner product of two matrices, that is $M\cdot N=\mbox{Tr}(MN^*)$ for $M,N\in\C^{d,d}$.
Let us apply now $\der$ to both members of the above equation: recall first our convention (\ref{cvpdt}) and the general relations 
$$
\der(\cj(dX^* \, m))=0,
\quad\mbox{ and }\quad
\der(\der X^* \, m)=-\der X^* \, \der m,
$$
which hold true for smooth processes $m$ and $X$.
Applying these relations to the right hand side of (\ref{exp1:irhodx}),
we end up with:
\begin{equation}\label{eq:di-rdx}
\der\lc\cj(dX^*\, r)\rc=\der X^* \, \der m - \der\xd \cdot \mu^* + \xd \cdot \der \mu^*
=\der X^* \, r + \xd \cdot \der \mu^*.
\end{equation}
When $m, X, \mu$ and $\xd$ are smooth enough, it is now clear that 
$\der[\cj(\der X^* \, r)]\in\cz\cac_3^{1+}$. Furthermore, this increment clearly satisfies the assumption (ii) of our  Proposition \ref{prop:Lambda}, and the fact that $\der[\cj(\der X^* \, r)]$ is an element of $\cac_3^{{\rm m},\beta}$ for any $\beta\le\ka$ easily stems from the second expression in (\ref{eq:di-rdx}). The increment $\der[\cj(\der X^* \, r)]$ is thus in the domain of application of $\laa$ (recall that $\der\der=0$), and we can write now 
$$
\cj(dX^* \, r)
=
\laa\lp \der X^* \, r + \xd \cdot \der \mu^* \rp.
$$
Reporting this identity into (\ref{exp2:imdx}), we end up with
\begin{equation}\label{exp3:imdx}
\cj(\der X^* \, m)
=
\der X^*\, m + \xd \cdot \mu^*
+\laa\lp \der X^* \, r + \xd \cdot \der \mu^* \rp.
\end{equation}

\vspace{0.3cm}

Observe that the expression above can be generalized to the non-smooth case,
since $\cj(\der X^* \, m)$ has been expressed now in terms of increments of $m$
and $X$. As a consequence, we shall use (\ref{exp3:imdx}) as a definition
of our extended integral, and summarize the previous considerations
in the following proposition:
\begin{proposition}\label{intg:mdx}
For a given $\ga>1/3$ and $1/3<\ka<\ga$,
let $X$ be a process satisfying Hypothesis \ref{hyp:x1}. Furthermore,  let
$m\in\cq_{\ka,b}(\Omega;\C^{d})$ such that $m_0=b$, $\delta m=\delta X^* \mu+r$,
 with $\mu\in\cac_1^\ka(\Omega;\C^{d,d})$ and $r\in\cac_2^{2\ka}(\Omega;\C^{d})$. Define $z$ by $z_0=a\in\C$ and
\begin{equation}\label{dcp:mdx}
\der z=
\der X^*\, m + \xd \cdot \mu^*
+\laa\lp \der X^* \, r + \xd \cdot \der \mu^* \rp.
\end{equation}
Finally, set
 \begin{equation*} \cj(d X^* \, m) = \der z. 
 \end{equation*}
Then:

\smallskip

\noindent
{\bf (1)}
$z$ is well-defined as an element of $\cq_{\ka,a}(\Omega;\C)$, and $\cj(d X^* \, m)$ coincides with the usual Riemann integral in case of two smooth processes $m$ and $X$. The quantity  $\cj_{ts}(d X^* \, m)$ is also defined as a Riemann integral for $t,s\in\Pi^+$.

\smallskip

\noindent
{\bf (2)}
Assume that $\oom\subset B(0,\tau)$, where $B(0,\tau)$ stands for the ball of radius $\tau$ centered at 0 in $\C$. Then the semi-norm of $z$ in $\cq_{\ka,a}(\Omega;\C)$ can be estimated as
\begin{equation}\label{bnd:norm-imdx-2}
\cn[z;\cq_{\ka,a}(\Omega;\C)]\le
c_{X}
\lp |b| + \tau^{\ga-\ka}\cn[m;\cq_{\ka,b}(\Omega;\C^{d})]\rp,
\end{equation}
for a positive constant $c_{X}$ depending only on $X$ and $\xd$.
Furthermore, the constant $c_{X}$ can be bounded as follows: 
$$
c_{X}\le c \lp \cn[X;\cac_1^{\ga}(\Omega;\C^d)]+\cn[\xd;\cac_2^{2\ga}(\Omega;\C^{d,d})] \rp,
$$ 
for a universal constant $c$. 

\smallskip

\noindent
{\bf (3)}
It holds
\begin{equation}\label{rsums:imdx}
\cj_{ts}(d X^* \, m)
=\lim_{|D_{ts}|\to 0}\sum_{i=0}^n
\lc (\der X^*)_{t_{i+1}, t_{i}} m_{t_{i}}
+ \xd_{t_{i+1}, t_{i}} \cdot \mu^*_{t_{i}} \rc
\end{equation} for any $s,t\in \oom$,
where the limit is taken over all partitions
$D_{ts} = \{s=t_0,\dots,t_n=t\}$
of $[s,t]$, as the mesh of the partition goes to zero.
\end{proposition}

\begin{proof}
we shall decompose this proof in two steps.

\vspace{0.3cm}

\noindent
{\it Step 1:}
Recalling the assumption $1/3<\ka<\ga$,
let us analyze the three terms in the right-hand side of (\ref{dcp:mdx})
and show that they define an element of $\cq_{\ka,a}$ such that
$\der z =\der X^*\, \zeta +r$ with
$$
\zeta=m
\quad\mbox{ and }\quad
\hat r= \xd \cdot \mu^*
+\laa\lp \der X^* \, r + \xd \cdot \der \mu^* \rp.
$$
Indeed, on the one hand, $m\in\cac_1^\ka(\Omega;\C^{d})$ and thus $\ze=m$ is of
the desired form for an element of $\cq_{\ka,a}$. On the other hand,
if $m\in\cq_{\ka,b}$, $\mu$ is assumed to be bounded and since
$\xd\in\cac_2^{2\ga}(\Omega;\C^{d,d})$ we get that
$\mu\cdot\xd\in\cac_2^{2\ga}(\Omega;\C)$.
Along the same lines we can prove that $r\,\der x\in \cac_3^{2\ka+\ga}(\Omega;\C)$
and $\xd \cdot \der\mu^*\in \cac_3^{\ka+2\ga}(\Omega;\C)$. Since
$\ka+2\ga\ge 2\ka+\ga>1$, we obtain  that
$\der X^* \, r + \xd \cdot \der \mu^*\in\cac_3^{2\ga+\ka}(\Omega;\C)$. In order to show that the latter increment is an element of $\dom(\laa)$, let us observe that the assumption (ii) of Proposition \ref{prop:Lambda} is easily satisfied, and we already mentioned that relation (\ref{eq:di-rdx}) entails $\der X^* \, r + \xd \cdot \der \mu^*\in\cac_3^{{\rm m},\beta}(\Omega;\C)$ for any $\beta\le\ka$. Hence, $\der X^* \, r + \xd \cdot \der \mu^*\in\dom(\laa)$, 
and
$$
\laa\lp \der X^* \, r + \xd \cdot \der \mu^* \rp\in \cac_2^{3\ka}(\Omega;\C).
$$
Thus we have proved that
$$
\hat r= \xd \cdot \mu^*+\laa\lp \der X^* \, r + \xd \cdot \der \mu^* \rp 
\in \cac_2^{2\ka}(\Omega;\C)
$$
and hence that $z\in\cq_{\ka,a}(\Omega;\C)$. The estimate (\ref{bnd:norm-imdx-2})
is now obtained using the same kind of considerations and is left
to the reader for the sake of conciseness. The analyticity $z,\zeta$ and $r$ is also a matter of standard considerations, as in Proposition \ref{cp:weak-phi}.

\vspace{0.3cm}

\noindent
{\it Step 2:}
The same kind of computations as those leading to (\ref{eq:di-rdx})
also show that
$$
\der\lp \der X^*\, m + \xd \cdot \mu^*  \rp
=
-\lc \der X^* \, r + \xd \cdot \der \mu^* \rc.
$$
Hence equation (\ref{dcp:mdx}) can also be read as
$$
\cj(dX^*\, m)=\lc \id-\laa\der \rc (\der X^*\, m + \xd \cdot \mu^* ),
$$
and a direct application of Proposition \ref{prop:integration} yields
(\ref{rsums:imdx}), which ends our proof.

\end{proof}

\begin{remark}
The previous proposition has a straightforward multidimensional
extension, which won't be stated here for sake of conciseness. For the nice continuity properties of the integral with respect to the driving process $X$, we refer to \cite[p. 101]{Gu}.
\end{remark}

\subsection{Rough differential equations of order 2}\label{sec:stoch-rough}

Recall that we wish to solve equations of the form (\ref{red:diff-eq}).
In our algebraic setting, we shall rephrase this as follows: we shall
say that $y$ is a local solution to (\ref{red:diff-eq}), if $y_0=a$,
$y\in\cq_{\ka,a}(\Omega;\C^l)$ and if there exists a neighborhood $\oom_0$ of 0 in $\bar\Pi^+$ such that, for any $s,t\in\oom_0$ we have
\begin{equation}\label{eds:alg-form}
(\der y)_{ts}=\cj_{ts}(dX^* \, \si(y)),
\end{equation}
where the integral $\cj(dX^*\, \si(y))$ has to be understood in the sense
of Proposition \ref{intg:mdx}. Our existence and uniqueness result
reads as follows:
\begin{theorem}\label{thm:ex-uniq}
Let $X$ be a process satisfying Hypothesis \ref{hyp:x1} and
$\si:\C^l\to\C^{d,l}$ be an analytic function. Then

\smallskip

\noindent
{\bf (1)}
There exists a neighborhood $\oom_0$ of 0 in $\bar\Pi^+$ such that
equation (\ref{eds:alg-form}) admits a unique solution $y$ in
$\cq_{\ka,a}(\oom_0;\C^l)$ for any $1/3<\ka<\ga$.

\smallskip

\noindent
{\bf (2)}
The mapping $(a,x,\xd)\mapsto y$ is continuous from
$\C^l\times\cac_1^{\ga}(\Omega_0;\C^d)\times\cac_2^{2\ga}(\Omega_0;\C^{d,d})$
to $\cq_{\ka,a}(\oom_0;\C^l)$.
\end{theorem}

\begin{proof}
We just sketch the proof of this theorem for sake of completeness:
we shall identify the solution on a small neighborhood $\oom_0$ as
the fixed point of the map $\Theta:\cq_{\ka,a}(\Omega_0;\C^l)\to\cq_{\ka,a}(\Omega_0;\C^l)$
defined by $\Theta(z)=\hz$ with $\hz_0=a$ and $\der\hz=\cj(dX^*\, \si(y))$.
The first step in this direction is to show that the ball
\begin{equation}\label{def:ball-m}
Q_M=\lcl z; \,z_0=a, \,\cn[z;\cq_{\ka,a}(\Omega_0;\C^l)]\le M \rcl
\end{equation}
is invariant under $\Theta$ for a certain $M>0$, if $\oom_0\subset B(0,\tau)$, with $\tau$ small enough. Indeed, due to  Propositions
\ref{cp:weak-phi} and \ref{intg:mdx} and assuming $\tau \leq 1$
we have
\begin{equation}\label{bnd:ggaz}
\cn[\Gamma(z);\cq_{\ka,a}(\Omega_0;\C^l)]\le c
\lp 1+c_{M}\tau^{\ga-\ka}(\cn[z;\cq_{\ka,a}(\Omega_0;\C^l)]+\cn^2[z;\cq_{\ka,a}(\Omega_0;\C^l)] )\rp.
\end{equation}
Since the set ${\mathcal A}=\{u\in\R_+^*:\,\,c( 1+c_{u}\tau^{\ga-\ka}(u+u^2))\le u\}$
is not empty as soon as $\tau$ is small enough, it is easily
shown that the ball $Q_M$ defined in (\ref{def:ball-m}) is left
invariant by $\Theta$ for $\tau$ small enough and $M$ in ${\mathcal A}$.

\smallskip

Now, since we are working in $Q_M$, the fixed point argument for $\Theta$ on $\oom_0$
is  standard  and  left to the reader (see \cite[Proposition 7]{Gu}). 

\end{proof}

\section{Rough path analysis: the general case}\label{sec:rough-general}
Having understood the rough path type tools we use for the case of a driving process $X$ with Hölder regularity $\ga>1/3$, we now proceed to a mere description of the generalization to the case of an arbitrary regularity $\ga>0$.
As in Lyons' theory, the algebraic integration setting for the resolution of rough differential equations relies on the a priori definition of a number
 of iterated integrals of the driving process $X$, which is shown to generate all the useful information needed to solve differential systems. More specifically, the following set of hypotheses is a generalization of the assumptions made in \cite{Gu2} to the complex plane:
\begin{hypothesis}\label{hyp:x}
Let $\oom$ be a neighborhood of 0 in $\Pi^+$, 
$X:\oom\to\C^d$ be a $\ga$-Hölder path with $0<\ga<1$, whose restriction to $\Pi^+$ is analytic, and set $N=\lfloor 1/\ga \rfloor$.
 Then $X$ is assumed to generate a family $\{\mathbf{X^n};\, 1\le n \le N\}$ defined on $\oom$, where $\mathbf{X^n}$ is a 1-increment with Hölder regularity $n\ga$, taking values in  $\C^{d^n}$, that is $\mathbf{X^n}\in\cac_2^{n\ga}(\Omega;\C^{d^n})$. We also suppose that $\mathbf{X^n}$ belongs to the space $\cac_2^{{\rm m},\ga}(\oom;\C^{d^n})$ defined by relation (\ref{eq:13}).
By definition $\mathbf{X^1}:=X$, and the algebraic relations satisfied by the $\mathbf{X^n}$ are the same as those of Proposition \ref{prop:dissec}: for any $n\le N$, $(i_1,\ldots,i_n)\in\{1,\ldots,d\}^n$, we have the multiplicative property
\beq\label{eq:chen-relation}
\der \mathbf{X^n}(i_1,\ldots,i_n)=\sum_{j=1}^{n-1} \mathbf{X^j}(i_1,\ldots,i_j) \, \mathbf{X^{n-j}}(i_{j+1},\ldots,i_n),
\eeq
that is $\der \mathbf{X^n}_{tus}(i_1,\ldots,i_n)=\sum_{j=1}^{n-1} \mathbf{X^j}_{tu}(i_1,\ldots,i_j) \, 
\mathbf{X^{n-j}}_{us}(i_{j+1},\ldots,i_n)$, for any $t,u,s\in\oom$. Furthermore, the rough path generated by $X$ is
said to be  {\em of geometric type} under the following additional condition: for any $n,m$ such that $n+m\le N$, we have:
\beq\label{eq:geom-rough-path}
\mathbf{X^n}(i_1,\ldots,i_n)\circ \mathbf{X^m}(j_1,\ldots,j_m)
=\sum_{\bar k\in\mbox{{\tiny Sh}}(\bar\imath,\bar\jmath)} \mathbf{X^{n+m}}(k_1,\ldots,k_{n+m}),
\eeq
where, for two tuples $\bar\imath,\bar\jmath$, $\Sigma_{(\bar\imath,\bar\jmath)}$ stands for the set of permutations of the indices contained in $(\bar\imath,\bar\jmath)$, and $\mbox{Sh}(\bar\imath,\bar\jmath)$ is a subset of $\Sigma_{(\bar\imath,\bar\jmath)}$ defined by:
$$
\mbox{Sh}(\bar\imath,\bar\jmath)=
\lcl  \si\in \Sigma_{(\bar\imath,\bar\jmath)}; \, 
\si \mbox{ does not change the orderings of } \bar\imath \mbox{ and } \bar\jmath \rcl.
$$
\end{hypothesis}
It should be mentioned at this point that relation (\ref{eq:chen-relation}) is just a version of Chen's relation for iterated integrals, while equation (\ref{eq:geom-rough-path}) states that any product of iterated integrals can also be expressed as a sum of iterated integrals.

\smallskip

For the sake of completeness, let us say a few words about the integration theory we rely on under Hypothesis \ref{hyp:x} (we refer to \cite{Gu2,TT} for further details). First of all, the class of processes we are able to integrate with respect to $X$ are called {\em analytic controlled processes of order $N-1$},
 and are defined as follows: we say that a path $z\in\cac_1^{\ka}(\oom;\C^{k})$, with $\frac{1}{N}<\ka<\ga<\frac{1}{N-1}$,
 is an {\em  analytic controlled path of order $N-1$}
 if its increments $(\der z(1),\ldots,\der z(k))$ can be decomposed into:
\beq\label{eq:dcp-weak-process-1}
\der z(i)=\sum_{n=1}^{N-1} \mathbf{X^n}(a_1,\ldots,a_n) \, \zeta^n(a_1,\ldots,a_n;i) + \rho^0(i),
\eeq
where the paths $\zeta^n(i)\in\cac_1^\ka(\Omega;\C^{k\times n})$ can be decomposed themselves as:
\beq\label{eq:dcp-weak-process-2}
\der \zeta^n(a_1,\ldots,a_n;i)
=\sum_{l=1}^{N-1-n} \mathbf{X^{l}}(b_1,\ldots,b_l) \, \zeta^{l+n}(b_1,\ldots,b_l,a_1,\ldots,a_n;i) + \rho^n(a_1,\ldots,a_n;i),
\eeq
and where the remainder terms $\rho^0,\ldots,\rho^{N-1}$ have the following regularity: for any $n\le N-1$, we have $\rho^n\in\cac_2^{(N-n)\ka}(\Omega;\C^{k})$. As in the case of order 2 described in the last section, we assume that the paths $z,\zeta^n,\rho^n$ are analytic on $\Pi^+$. Denote by $\cq_{\ka,a}(\Omega;\C^{k})$, or simply  by $\cq_{\ka,a}$ for notational sake, the set of $\C^{k}$-valued analytic controlled paths. Then a natural norm on this space is given by:
$$
\cn[z;\cq_{\ka,a}]=\sum_{n=0}^{N-1} \cn[\zeta^n;\cac_1^\ka] 
+ \sum_{n=1}^{N-1} \cn[\zeta^n;\cac_1^0]
+\sum_{n=0}^{N-1} \cn[\rho^n;\cac_2^{(N-n)\ka}].
$$

\smallskip

With this notions in hand, it is worth recalling how to integrate controlled processes. This is summarized in the following proposition, which is a simplified version of \cite[Theorem 8.5]{Gu2} adapted to the complex plane, and for which we use the traditional convention of summation over repeated indices.
\begin{proposition}\label{prop:intg-z-dx}
For a given $\ga>0$, let $X$ be a process satisfying Hypothesis \ref{hyp:x}. Furthermore, for $\lceil 1/\ga\rceil^{-1}<\ka<\ga$, let $z\in\cq_{\ka,b}(\Omega;\C^{d})$, with a decomposition given by (\ref{eq:dcp-weak-process-1}) and (\ref{eq:dcp-weak-process-2}). Define $\hz$ by $\hz_0=a\in\C$ and
\beq\label{eq:dcp-intg-z-dx}
\der \hz=
\der x(i)\, z(i)+ \sum_{n=1}^{N-1} \mathbf{X^{n+1}}(a_1,\ldots,a_n,i) \, \zeta^n(a_1,\ldots,a_n;i) 
+\laa(U),
\eeq
with
$$
U=
\sum_{n=0}^{N-1} \mathbf{X^{n+1}}(a_1,\ldots,a_n,i)  \rho^n(a_1,\ldots,a_n;i)  
+ \mathbf{X^{N}}(a_1,\ldots,a_{N-1},i) \der\zeta^{N-1}(a_1,\ldots,a_{N-1};i).
$$
Finally, set $\cj(m\, dx) = \der \hz$. Then:

\smallskip

\noindent
{\bf (1)}
$\hz$ is well-defined as an element of $\cq_{\ka,a}(\Omega;\C)$, and $\cj(m\, dx)$ coincides with the usual Riemann integral in case of two analytic processes $m$ and $X$.

\smallskip

\noindent
{\bf (2)}
Assume that $\oom\subset B(0,\tau)$, where $B(0,\tau)$ stands for the ball of radius $\tau$ centered at 0 in $\C$.
Then the semi-norm of $\hz$ in $\cq_{\ka,a}(\Omega;\C)$ can be estimated as
\begin{equation}\label{bnd:norm-imdx}
\cn[\hz;\cq_{\ka,a}(\Omega;\C)]\le
c_{X}
\lp 1 + \tau^{\ga-\ka}\cn[z;\cq_{\ka,b}(\Omega;\C^{d})]\rp,
\end{equation}
for a positive constant $c_{X}$ depending only on $X^1,\ldots,X^N$. 

\smallskip

\noindent
{\bf (3)}
It holds
\begin{equation*}
\cj_{ts}(z\, dx)
=\lim_{|D_{ts}|\to 0}\sum_{k=0}^n
\lc 
\der X_{t_{k+1},t_{k}}(i)\, z_{t_{k}}(i)
+ \sum_{n=1}^{N-1} \mathbf{X^{n+1}}_{t_{k+1},t_{k}}(a_1,\ldots,a_n,i) \, 
\zeta^n_{t_{k}}(a_1,\ldots,a_n;i) 
 \rc
\end{equation*} 
for any $s,t\in\oom$, where the limit is taken over all partitions
$D_{ts} = \{s=t_0,\dots,t_n=t\}$
of $[s,t]$, as the mesh of the partition goes to zero.
\end{proposition}

\smallskip

Furthermore, the construction above allows to solve rough differential equations in a reasonable sense:
\begin{proposition}
Under the same conditions as for Proposition \ref{prop:intg-z-dx}, let $\si$ be a function from $\C^{n}$ to $\C^{n\times d}$, analytic in each of its variables. Then there exists a neighborhood $\oom_0$ of 0 in $\bar\Pi^+$ such that the differential equation
$$
(\der y)_{ts}= \cj_{ts}(\si(y)\, dx), \quad s,t\in\oom,
$$
where the integral in the right hand side has to be understood as in equation (\ref{eq:dcp-intg-z-dx}), admits a unique solution on $\oom_0$, living in the class of $\C^{n}$-valued analytic controlled processes.
\end{proposition}


\section{Analytic fractional Brownian motion and preliminaries}
We review in this section the construction of the analytic fBm $\gga$, and we also include here some useful preliminary results concerning the regularity of increments in $\oom$, and some complex analysis estimates for the kernel $(-\imath(z-\bar w))^{2\al-2}$ defined on $\Pi^+$.
\subsection{Definition of the analytic fBm}
\label{sec:def-fBm}


As mentioned in the Introduction, the article \cite{Un} is an elaboration of 
a stochastic calculus with respect to the fractional Brownian motion by analytic continuation.
 More specifically, a complex-valued processed indexed by $z\in\Pi^+$, called $\gga$, is introduced there.
 This process is analytic on $\Pi^+$ and converges in every reasonable sense to a continuous
process with real-time parameter (still denoted by $\Gamma$) when the imaginary part
of $z$ goes to $0$. The current section is devoted to recall this formalism,
 which is the one we shall adopt in order to construct a fractional rough path for $\Gamma$
 for any Hurst parameter $\al\in(0,1)$.
 Notice that, since the case $\al>1/2$ is trivial from the rough path analysis point of view,
 we shall  assume in the sequel that  $\al\in(0,1/2)$. The Brownian case $\alpha=1/2$ may be seen as a limit.

\smallskip

Let us first recall some classical notations of complex analysis: for $x\in\R$ and $k\in\N$, the Pochhammer symbol $(x)_k$ is defined by:
$$
(x)_k=\prod_{j=0}^{k-1} (x+j)=\frac{\mathbf{\gga}(x+k)}{\mathbf{\gga}(x)},
$$ 
where $\mathbf{\gga}$ stands for the usual Gamma function. Recall that we denote by
 $\Pi^+=\{ z=x+\imath y\ |\ x\in\R, y>0\}$, resp. $\bar{\Pi}^+=\{z=x+\imath y\ |\ x\in\R, y\ge 0\}$
the open, resp. closed  upper half-plane in $\C$. Similarly, $\Pi^-$, resp. $\bar{\Pi}^-$ stand
for the open, resp. closed lower half-planes.

\smallskip

With these notations in mind, the easiest way to define $\gga=\{\gga_z;\, z\in\Pi^+\}$ makes use of a series expansion involving the analytic functions $\{f_k;\, k\ge 0\}$, defined on $\pp$ by:
$$
f_k(z)=2^{\al-1} \lc  \frac{\al(1-2\al)(2-2\al)_k}{2\cos(\pi\al) k!} \rc^{1/2}
\lc  \frac{z+\imath}{2\imath} \rc^{2\al-2}
\lc  \frac{z-\imath}{z+\imath} \rc^{k}.
$$
It is shown in \cite{Un} that the series $\sum_{k\ge 0}f_k(z) \overline{f_k(w)}$ converges in absolute value for $z,w\in\pp$, and that the following identity holds true:
\beq\label{eq:dcp-cov-gamma}
\sum_{k\ge 0}f_k(z) \overline{f_k(w)} 
= \frac{\al(1-2\al)}{2\cos(\pi\al)} \lp -\imath (z-\bar w) \rp^{2\al-2}.
\eeq
This fact allows to define the process $\Gamma$ in the following way:
\begin{proposition}
Let $\{\xi_k^1,\xi_k^2;\, k\ge 0\}$ be two families of independent standard Gaussian random variables, defined on a complete probability space $(\cu,\cf,\bp)$, and for $k\ge 0$, set $\xi_k=\xi_k^1+\imath\xi_k^2$. Consider the process $\gga'$ defined for $z\in\pp$ by $\gga'_z=\sum_{k\ge 0}f_k(z)\xi_k$. Then:

\smallskip

\noindent
{\bf (1)}
$\gga'$ is a well-defined analytic process on $\pp$.

\smallskip

\noindent
{\bf (2)}
Let $\gga:(0,1)\to\pp$ be any continuous path with endpoints $\gga(0)=0$ and $\gga(1)=z$, and set $\gga_z=\int_{\gga}\gga'_u \, du$. Then $\gga$ is an analytic process on $\pp$. Furthermore, as $z$ runs along any path in $\pp$ going to $t\in\R$, the random variables $\gga_z$ converge almost surely to a random variable called again $\gga_t$.

\smallskip

\noindent
{\bf (3)}
The family $\{\gga_t;\, t\in\R\}$ defines a Gaussian centered complex-valued process, whose covariance function is given by:
$$ \be[\Gamma_s \Gamma_t]=0, \quad 
\be[\gga_s\bar\gga_t]=
\frac{e^{-\imath\pi\al \,\sgn(s)}|s|^{2\al}
+e^{\imath\pi\al \,\sgn(t)}|t|^{2\al}-e^{\imath\pi\al \,\sgn(t-s)}|s-t|^{2\al}}{4\cos(\pi\al)}.
$$
The paths of this process are almost surely $\ka$-Hölder for any $\ka<\al$.

\smallskip

\noindent
{\bf (4)}
Both real and imaginary parts of $\{\gga_t;\, t\in\R\}$ are (non independent) fractional Brownian motions indexed by $\R$, with covariance given by
\begin{equation} \be [\Re \Gamma_s \Im \Gamma_t]=-\frac{\tan\pi\alpha}{8} \left[ -\sgn(s) |s|^{2\alpha}+
\sgn(t) |t|^{2\alpha}-\sgn(t-s) |t-s|^{2\alpha} \right].\end{equation}
\end{proposition}

\begin{remark}\label{rmk:def-B}
It should be stressed at this point that the paper \cite{Un} mainly focuses on the real part of $\gga$, that is a standard fractional Brownian motion. We shall see however that $\gga$ is an interesting process in its own right, insofar as it allows the construction of a rough path for any value of the Hurst parameter $\al\in(0,1)$.
\end{remark}

\smallskip

Let us also recall some basic facts about $\gga$ which will be used extensively in the sequel: first, according to (\ref{eq:dcp-cov-gamma}), the (Hermitian) covariance between $\gga'_z$ and $\gga'_w$ for $z,w\in\Pi^+$ is given by:
\beq\label{eq:cov-gamma-prime}
\be\lc \gga'_z \, \bar \gga'_w \rc = \frac{\al(1-2\al)}{2\cos(\pi\al)} \lp -\imath (z-\bar w) \rp^{2\al-2}.
\eeq
The following formula will be used throughout the article: for a piecewise smooth path $\gamma:(0,1)\to\pp$, we have:
\beq\label{eq:cov-gamma-prime-path}
\be\lc \int_\gamma \gga'_z dz \, \int_\gamma \bar \gga'_w dw \rc
= \frac{\al(1-2\al)}{2\cos(\pi\al)} \int_\gamma dz \int_\gamma dw \lp -\imath (z-\bar w) \rp^{2\al-2}.
\eeq

\subsection{Garsia-Rodemich-Rumsey type lemmas}\label{sec:grr}
This section is devoted to recall or give some deterministic regularity results for increments, which will be essential in order to quantify the convergence of the approximations of our process $\Gamma$. First let us recall a particular case of a classical lemma due to Garsia \cite[Lemma 2]{Ga}:
\begin{lemma}\label{lem:4.3}
Let $f$ be a continuous function defined on a compact set $D\subset\R^d$ for $d\ge 1$, and set, for $p\ge 1$
$$
U_{\ka,p}(f)=\lp \int_D\int_D \frac{|(\der f)_{wv}|^{2p}}{|w-v|^{2\ka p+2d}}\rp^{1/2p}.
$$
Then $\cn[f;\, \cac_1^{\ka}(D)]\le c \, U_{\ka,p}(f)$, for a universal positive constant $c$.
\end{lemma}

When $D=\oom\subset\bar{\Pi}^+$, we need an extension of this lemma to increments which are not necessarily written as $\der f$ for functions $f\in\cac_1$:
\begin{proposition}\label{prop:4.4}
Let $\Omega:=B(0,r)\cap \bar{\Pi}^+$ be a neighborhood of $0$ in $\bar{\Pi}^+$, and
${\mathcal R}\in C_2(\Omega;\C^n)$ for $n\ge 1$ such that $\delta {\mathcal R}\in  C_3^{\kappa}(\Omega;\C^n)$. Set for $p\ge 1$
\begin{equation}
U_{\ka,p}(\Omega;{\mathcal R}):=\lp\int_{\Omega}\int_{\Omega} \frac{|{\mathcal R}_{wv}|^{2p}}{|w-v|^{2\kappa p+4}}
\ dv \ dw\rp^{1/2p},
\end{equation}
and assume $U_{\ka,p}(\Omega;{\mathcal R})<\infty$. Then ${\mathcal R}\in C_2^{\kappa}(\Omega;\C^n)$; more precisely,
\begin{equation}
\cn[{\mathcal R};\, \cac_2^{\ka}(\Omega;\C^n)] \le 
c \lp U_{\ka,p}(\Omega;{\mathcal R})+\cn[\delta {\mathcal R};\, \cac_3^{\ka}(\oom;\C^n)]\rp,
\end{equation}
for a universal constant $c>0$.
\end{proposition}

\begin{proof}

Let $s,t\in \Omega':=B(0,r/4)\cap \bar{\Pi}^+\subset \Omega$. We wish to show that 
\begin{equation}\label{eq:50}
\crr_{ts}\le c \lp U_{\ka,p}(\Omega;{\mathcal R})+\cn[\delta {\mathcal R};\, \cac_3^{\ka}(\oom;\C^n)]\rp |t-s|^{\ka}.
\end{equation} 
To this end, let us construct a sequence of
points $(s_k)_{k\ge 0}$, $s_k\in\Omega$ converging to $t$ in the following way: set  $s_0=t$, suppose by induction that
$s_0,\ldots, s_k$ have been constructed, and let $V_k:=B(s,\frac{|s_k-s|}{2})\cap\bar{\Pi}^+$. Note that, since  we are working on the upper half plane $\Pi^+$, 
the area $\mu(V_k)$ of $V_k$ is at least $\half \mu(B(s,\frac{|s_k-s|}{2})=\frac{\pi}{8} |s_k-s|^2$. Define then
\begin{equation}
A_k:=\left\{v\in V_k\ |\ I_v >\frac{16}{\pi} \frac{U_{\ka,p}^{2p}(\Omega;{\mathcal R})}{|s_k-s|^2} \right\} \label{eq:Ak}
\end{equation}
and
\begin{equation}
B_k:=\left\{v\in V_k\ |\ \frac{|{\mathcal R}_{s_k v}|^{2p}}{|s_k-v|^{2\kappa p+4}}>\frac{16}{\pi |s_k-s|^2}
I(s_k) \right\} \label{eq:Bk}
\end{equation}
where we have set
\begin{equation*} I(v):=\int_{B(s,|v-s|)\cap\bar{\Pi}^+} \frac{|{\mathcal R}_{uv}|^{2p}}{|v-u|^{2\kappa p+4}}
\ du.
\end{equation*}
Let us prove now that $V_k\setminus (A_k\cup B_k)$ is not empty: observe that
\begin{equation*}
U_{\ka,p}^{2p}(\Omega;{\mathcal R}) \ge \int_{A_k} dv I(v) >\frac{16}{\pi} \frac{U_{\ka,p}^{2p}(\Omega;
{\mathcal R})}{|s_k-s|^2} \mu(A_k) 
\end{equation*}
and
\begin{equation*}
I(s_k)\ge \int_{B_k} \frac{|{\mathcal R}_{us_k}|^{2p}}{|s_k-u|^{2\kappa p+4}} \ du>\frac{16}{\pi}
\frac{\mu(B_k)}{|s_k-s|^2} I(s_k).
\end{equation*}
All together one has obtained $\mu(A_k),\mu(B_k)<\frac{\pi}{16} |s_k-s|^2$ so $\mu(A_k)+\mu(B_k)<\mu(V_k)$.
One now chooses $s_{k+1}$ arbitrarily in $V_k\setminus (A_k\cup B_k)$. Note that, by construction,
$|t-s|<r/2$, and $\left| \frac{s_{k+1}-s}{s_k-s} \right|<1/2$ so $s_k\to s$ while staying inside $\Omega$.

\smallskip

Now decompose (by using a number of times the operator $\delta$) ${\mathcal R}_{s_0,s}$ into
\begin{equation} {\mathcal R}_{s_0 s}={\mathcal R}_{s_{n+1} s}+\sum_{k=0}^n \left(
{\mathcal R}_{s_k s_{k+1}} +(\delta {\mathcal R})_{s_k s_{k+1} s} \right). \label{eq:dec-del} \end{equation}
Applying (\ref{eq:Bk})$_k$ and (\ref{eq:Ak})$_{k-1}$, one gets
\begin{equation*} \frac{ |{\mathcal R}_{s_k s_{k+1}}|^{2p}}{|s_k-s_{k+1}|^{2\kappa p+4}} <
\frac{16}{\pi |s_k-s|^2} \ .\ \frac{16}{\pi} \frac{ U_{\ka,p}^{2p}(\Omega;{\mathcal R})  }{|s_{k-1}-s|^2}
<\frac{256}{\pi^2}  U_{\ka,p}^{2p}(\Omega;{\mathcal R})  |s_k-s|^{-4}.
\end{equation*}
Recalling our convention $a\lesssim b$ for the relation $a\le C \,b$, where $C$ is a given universal constant, we obtain $|{\mathcal R}_{s_k s_{k+1}}|\lesssim U_{\ka,p}(\Omega;\crr)  |s_k-s|^{\kappa}$. Furthermore, we have by construction $|s_k-s|\lesssim 2^{-n}|t-s|$, and thus
\begin{equation} \label{eq:54}
\left| \sum_{k=0}^n {\mathcal R}_{s_k, s_{k+1}}\right|\lesssim U_{\ka,p}(\Omega;\crr) |t-s|^{\kappa}. \end{equation}

\smallskip

Turning now to $\delta {\mathcal R}$, it is easily seen that $|\delta {\mathcal R}_{s_k s_{k+1},s}|\lesssim \cn[\der\crr;\, \cac_3^{\ka}(\oom;\C^n)] |s_k-s|^{\kappa}$. Invoking again the relation $|s_k-s|\lesssim 2^{-n}|t-s|$, we end up with
\begin{equation} \label{eq:55}
\sum_{k=0}^n  |\delta {\mathcal R}_{s_k s_{k+1},s}|\lesssim  
\cn[\der\crr;\, \cac_3^{\ka}(\oom;\C^n)] \, |t-s|^{\kappa}.
\end{equation}

\smallskip

Finally, plugging relations (\ref{eq:54})-(\ref{eq:55}) into (\ref{eq:dec-del}) and  letting $n\to \infty$, we easily get the announced bound  (\ref{eq:50}), which ends the proof.

\end{proof}

\subsection{Complex analysis preliminaries}
\label{sec:complex-analysis}

The identity (\ref{eq:cov-gamma-prime-path}) involves integrals along  some piecewise smooth paths in $\C$, which have to be estimated. We summarize in this section the upper bounds which will be needed later on. 

\smallskip

First of all, the integral appearing in (\ref{eq:cov-gamma-prime-path}) can be estimated thanks 
to the following lemma borrowed from \cite[Lemma 1.5]{Un}:
\begin{lemma}\label{lem:bound-intg-path}
Let $\ga:(0,1)\to\pp$ be a piecewise smooth, continuous path. Then
$$
\lln \int_\ga dz \int_{\bar{\ga}} d\bar{w} \lp -\imath (z-\bar w) \rp^{2\al-2} \rrn
\le c |\ga(1)-\ga(0)|^{2\al},
$$
for a universal positive constant $c$.
\end{lemma}

The following bound on iterated integrals of the process $\gga$, shown in \cite[Theorem 3.4]{Un} (where they
are called {\em analytic iterated integrals})
 can then be  seen as an extension of the previous lemma. 
\begin{lemma}[\textbf{analytic iterated integrals}] \label{lem:cvg-anal-intg}
Consider $s,t$ in a fixed bounded  neighborhood of $0$ in $\bar{\Pi}^+$, and let $f_1,\ldots,f_n$ and $g_1,\ldots,g_n$ be analytic functions defined on a neighborhood $V$ of the closed strip $\bar{\Pi}_{s,t}^+:=\{z\in\C\ |\ z=\lambda s+(1-\lambda)t+\imath\mu|t-s|, \ \
\lambda,\mu\in[0,1]\}$. Let also $\gga=(\gga(1),\ldots,\gga(d))$ be a $d$-dimensional analytic
fractional Brownian motion, where each component $\gga(j)$ is defined as in Section \ref{sec:def-fBm}. For $\ep,\eta$ small enough, define $\cv_{s,t}(\ep,\eta)$ by
\begin{align*}
&\cv_{ts}(\ep,\eta)  =\be \left[ Z_1 \, Z_2\right]  \\
&= \left( \int_{[s;t]} du_1 \int_{[s,u_1]} du_2\ldots \int_{[s,u_{n-1}]} du_n\right)
 \left( \int_s^t dv_1 \int_{[s,v_1]} dv_2\ldots \int_{[s,v_{n-1}]} dv_n\right) \\
&\hspace{5cm}
\times \prod_{j=1}^{n}  f_j(u_j+\im \ep) \, \overline{ g_j(v_j+\im\eta)}
\, (-\im (u_j-\bar{v}_j)+\ep+\eta)^{2\al-2}\, du_j\, dv_j,
\end{align*}
where $Z_1$ is defined by
$$
\int_{[s;t]}  f_1(u_1+\imath\eps) d\Gamma_{u_1+\imath\eps}(1) \int_{[s;u_1]} f_2(u_2+\imath\eps) d\Gamma_{u_2+\imath\eps}(2) \ldots
\int_{[s;u_{n-1}]} f_n(u_n+\imath\eps) d\Gamma_{u_n+\imath\eps}(n),
$$
and $Z_2$ can be written as:
$$
\int_{[s;t]}  \overline{g_1(v_1+\imath\eps)} d\bar{\Gamma}_{v_1+\imath\eps}(1) \int_{[s;v_1]} \overline{g_2(v_2+\imath\eps)}
 d\bar{\Gamma}_{v_2+\imath\eps}(2)\ldots
\int_{[s;v_{n-1}]} \overline{g_n(v_n+\imath\eps)} d\bar{\Gamma}_{v_n+\imath\eps}(n).
$$
Then the following bound holds true:
$$
|\cv_{ts}(\ep,\eta)| \lesssim \prod_{j=1}^{n} \sup_{z\in \bar{\Pi}^+_{s,t}} |f_j(z)| \,
\prod_{j=1}^{n} \sup_{z\in \bar{\Pi}^+_{s,t}} |g_j(z)|
\, |t-s|^{2\al n}.
$$
\end{lemma}

The last ingredient we need for our computations is a specific bound for analytic functions 
integrated with respect to the kernel $(-\im(x-y))^{2\al-2}$. Observe that this bound will 
not be used directly in the sequel, but will serve as a prototype for our future computations.

\begin{lemma}\label{lem:bound-intg-contour}
Let $s,t$ in a fixed bounded neighborhood of $0$ in $\bar{\Pi}^+$, and $\phi(z,\bar{w})$ be an analytic function on a neighborhood of
 $\bar{\Pi}^+_{(s,t)}\times \bar{\Pi}^-_{(s,t)}$, where  
\begin{equation} \bar{\Pi}_{s,t}^+:=\{z\in\C\ |\ z=\lambda s+(1-\lambda)t+\imath\mu|t-s|, \ \
\lambda,\mu\in[0,1]\}, \quad \bar{\Pi}_{s,t}^-:=\{\bar{z}\ |\ z\in\bar{\Pi}_{s,t}^+\} \end{equation}

For $\eps,\eta>0$, define $\Theta(\phi)$ as:
\begin{multline*}
\lc \Theta(\phi)\rc (\eps,\eta;s,t)\\
:=\int_{[s,t]} dz \int_{[\bar s,\bar t]} d\bar{w} \left[ (-\imath(z-\bar{w})+2\eps)^{2\alpha-2} - (-\imath(z-\bar{w})+\eps+\eta)^{2\alpha-2}
 \right] \phi(z,\bar{w}) 
\end{multline*}
Then, for every $\rho\in(0,2\alpha)$, there exists $C_{\rho}$ such that
\begin{equation}
 |\lc \Theta(\phi)\rc(\eps,\eta;s,t)|\le 
 C_{\rho} |\eps-\eta|^{\rho} |t-s|^{2\alpha-\rho} \, M^\phi_{ts},
 \label{eq:bound} 
 \end{equation}
 where
 $$
 M^\phi_{ts}\triangleq  
 \sup\lcl  |\phi(u,\bar{v})|; \, (u,\bar{v})\in\bar{\Pi}^+_{(s,t)}\times \bar{\Pi}^-_{(s,t)}\rcl.
 $$
\end{lemma}

\begin{proof}
Without restriction of generality we may assume that $\eps>\eta>0$. We use the following contour of integration in
$\bar{\Pi}^+_{(s,t)}\times \bar{\Pi}^-_{(s,t)}$:
\begin{equation} \Delta:=\Gamma\times\bar{\Gamma},\quad \Gamma=\Gamma_1\cup\Gamma_2\cup \Gamma_3:=[s,s+\imath|t-s|]\cup[s+\imath|t-s|,t+\imath|t-s|]\cup
[t+\imath|t-s|,t].\end{equation}
Set $\Delta_{i,j}=\Gamma_i\times \bar{\Gamma}_j$ so that $\Delta=\cup_{1\le i,j\le 3} \Delta_{i,j}$. Let $I_{i,j}$ be the integral over
$\Delta_{i,j}$ of the function $(z,\bar{w})\mapsto   \left[ (-\imath(z-\bar{w})+2\eps)^{2\alpha-2} - 
(-\imath(z-\bar{w})+\eps+\eta)^{2\alpha-2} \right] \phi(z,\bar{w}) $. We shall give a bound of type  (\ref{eq:bound}) for each $I_{i,j}$. The proof relies on the following observation: if $\eps,\eta>0$ and $z\in \C$, $\Re z>0$,   then
(for any $\rho\in(0,1)$)
\begin{equation} |(z+\eps)^{2\alpha-2}-(z+\eta)^{2\alpha-2}|\le C|\eps-\eta|^{\rho} |z|^{2\alpha-2-\rho}.\end{equation}

 By symmetry we only need to consider
the following four cases (only the fourth one is non-trivial since $z$ and $\bar{w}$ may be $\eps$-close) :

\smallskip

\noindent
{\it Case 1:} $i=j=2$.
\begin{eqnarray}\label{eq:bnd-I22}
 |I_{2,2}|&=&\Big| \int_{[s;t]} dz \int_{[\bar s;\bar t]} d\bar{w} \big[ (-\imath(z-\bar{w})+2|t-s|+2\eps)^{2\alpha-2} 
\nonumber\\
&& \hspace{2cm} -(-\imath(z-\bar{w})+2|t-s|+\eps+\eta)^{2\alpha-2} \big]
\phi(z+\imath|t-s|,\bar{w}-\imath|t-s|) \Big| \nonumber\\
&\le& C_{\rho} \int_{[s;t]} |dz| \int_{[\bar s; \bar t]} |d\bar{w}| |t-s|^{2\alpha-2-\rho} (\eps-\eta)^{\rho} \, M^\phi_{ts}
=C'_{\rho} |t-s|^{2\alpha-\rho} (\eps-\eta)^{\rho} \, M^\phi_{ts}.
\end{eqnarray}

\smallskip

\noindent
{\it Case 2:} $i=1,j=3$.
\begin{eqnarray}\label{eq:bnd-I13}
 |I_{1,3}|&=&\Big|  \int_0^{|t-s|} dx \int_0^{|t-s|} dy \big[ (-\imath(s-\bar{t})+x+y+2\eps)^{2\alpha-2} 
 \nonumber\\
&& \hspace{4cm}  -(-\imath(s-\bar{t})+x+y+\eps+\eta)^{2\alpha-2}\big]
\phi(s+\imath x,t-\imath y) \Big| \nonumber\\
& \le& C_{\rho} \int_0^{|t-s|} dx \int_0^{|t-s|} dy |t-\bar{s}|^{2\alpha-2-\rho} (\eps-\eta)^{\rho} \, M^\phi_{ts}
\nonumber\\
&\le &  C'_{\rho} |t-s|^{2\alpha-\rho} (\eps-\eta)^{\rho} \, M^\phi_{ts},
\end{eqnarray}
since $|t-\bar{s}|\ge |t-s|$.

\smallskip

\noindent
{\it Case 3:} $i=1,j=2$.
\begin{eqnarray}\label{eq:bnd-I12}
 |I_{1,2}|&=&\Big| \int_{0}^{|t-s|} dx \int_{[\bar s; \bar t]} d\bar{w} \big[ (-\imath(s-\bar{w})+|t-s|+x+2\eps)^{2\alpha-2} 
 \nonumber\\
&& \hspace{2cm} - (-\imath(s-\bar{w})+|t-s|+x+\eps+\eta)^{2\alpha-2}
\big] \phi(s+\imath x,\bar{w}-\imath|t-s|) \Big| \nonumber\\
& \le& C_{\rho} \int_0^{|t-s|} dx \int_{[\bar s; \bar t]} |d\bar{w}| |t-s|^{2\alpha-2-\rho} (\eps-\eta)^{\rho} \, M^\phi_{ts}=C'_{\rho} |t-s|^{2\alpha-\rho} (\eps-\eta)^{\rho} \, M^\phi_{ts}. \nonumber\\
\end{eqnarray}

\smallskip

\noindent
{\it Case 4:} $i=1,j=1$.
\begin{eqnarray}\label{eq:bnd-I11}
 |I_{1,1}|&=&\left| \int_0^{|t-s|} dx \int_0^{|t-s|} dy \left[ K_{2\eps}(x,y)- K_{\eps+\eta}(x,y) \right]
\phi(s+\imath x,s-\imath y) \right|\\
& \le& C_{\rho} \int_0^{|t-s|} dx \int_0^{|t-s|} dy (x+y)^{2\alpha-2-\rho} 
(\eps-\eta)^{\rho} \, M^\phi_{ts} \le C'_{\rho} |t-s|^{2\alpha-\rho} (\eps-\eta)^{\rho} \, M^\phi_{ts}, 
\nonumber
\end{eqnarray}
where we have set $K_{a}(x,y)=(2\Im s+x+y+a)^{2\alpha-2}$ for any positive $a$.
It should be noticed at this point that the last integral converges only if $\rho<2\alpha$, which is one of our standing assumptions. Now, putting together the estimates (\ref{eq:bnd-I22}), (\ref{eq:bnd-I13}), (\ref{eq:bnd-I12}) and (\ref{eq:bnd-I11}), we get the desired result.

\end{proof}

\begin{remark} The kernel $(x,y)\to (x+y)^{2\alpha-2-\rho}$ appearing in the last case is singular only
at the point $(x,y)=(0,0)$,
whereas the usual kernel $(x,y)\to (\pm\imath(x-y)))^{2\alpha-2}$ is singular on the diagonal. This simple fact explains why
our estimates work (and why the deformation of contour is so important). Note that the absolute
value should {\em not} be placed inside the integral {\em before} the deformation of contour
(otherwise the integrals become most of the time  infinite in the  limit $\eta\to 0$).

\end{remark}

\section{The rough path associated to $\Gamma$}
We proceed in this section to the definition of a rough path above the process $\gga$ defined at Section \ref{sec:def-fBm}. As mentioned in the introduction, this will be achieved by regularizing $\gga$ into a process $\gga^\ep$ defined on $\Pi^+$ by $\gga^\ep_t=\gga_{t+\imath \ep}$. This latter process is analytic on $\Pi^+$, which allows to define any iterated integral of $\gga^\ep$ in the Riemann sense. Then the convergence of these integrals in some suitable Hölder spaces is obtained by combining the Garsia type result of Proposition \ref{prop:4.4} and some moment estimates similar to Lemma \ref{lem:bound-intg-contour}.

\subsection{Convergence of $\mathbf{\gga^{\ep}}$}
\label{sec:cvgce-B-ep}

A very first step in the analysis of $\gga$ consists in getting some convergence results for $\gga^{\ep}$ itself towards $\gga$, in Hölder spaces. In order to obtain this (intuitively trivial) convergence, we shall use the following elementary estimate:

\begin{lemma}\label{lem:moments-B-ep}
For all $\ep,\eta >0$, and $s,t \in\Pi^+$, we have
$$
\be\lc \lln \gga^{\ep}_t-\gga^{\ep}_s\rrn^2 \rc \le c \lln t-s\rrn^{2\al} ,
\quad\text{and}\quad
\be\lc \lln \gga^{\ep}_t-\gga^{\eta}_t\rrn^2 \rc \le  c \lln \ep-\eta \rrn^{2\al},
$$
where the constants $c$ do not depend on $\ep,\eta,s,t$.
\end{lemma}

\begin{proof}
For the first inequality, observe that
$$
\gga^{\ep}_t-\gga^{\ep}_s = \int_{[s+\imath\eta,t+\imath\eta]} \gga'_z \, dz,
$$
so that
$$
E\lc \lln \gga^{\ep}_t-\gga^{\ep}_s \rrn^2 \rc  \leq 
c \lln \int_{[s+ \imath\eta,t+ \imath\eta]} \int_{[\bar s-\imath\eta,\bar t-\imath\eta]} 
(-\imath(z-\bar w))^{2\al-2} dz d\bar{w} \rrn 
\leq  c \lln t-s \rrn^{2\al},
$$
owing to Lemma \ref{lem:bound-intg-path}. For the second inequality, use the decomposition:
$$
\gga^{\ep}_t-\gga^{\eta}_t =\int_{[t+ \imath\ep,t+ \imath\eta]} \gga'_z \, dz,
$$
which yields, with the same kind of arguments,
$$
E\lc \lln \gga^{\ep}_t-\gga^{\eta}_t \rrn^2 \rc  \leq  c 
\lln \int_{[t+ \imath\ep,t+ \imath\eta]} \int_{[\bar t-\imath\ep,\bar t-\imath\eta]} (-\imath(z-\bar w))^{2\al-2}
dz d\bar{w}  \rrn 
\leq  c \lln \ep-\eta \rrn^{2\al}.
$$

\end{proof}

We are now ready to study the convergence of $\gga^{\ep}$ on our fixed neighborhood $\oom$ (recall also that we work on a complete probability space $(\cu,\cf,\bp)$):
\begin{lemma}\label{lemma:cvgce-B-ep}
As $\ep\to 0$, the process $\gga^{\ep}$ converges in $L^1(\cu;\cac_1^\ga(\oom))$, for any $\ga<\al$ and $T>0$. Its limit is the analytic fractional Brownian motion $\gga$.
\end{lemma}

\begin{proof}
We shall divide this proof into two steps:

\smallskip

\noindent
{\it Step 1: Reduction to moment estimates.}
We shall prove that $\{\gga^{\ep};\,\ep>0\}$ is a Cauchy sequence in $L^1(\cu;\cac_1^\ga(\oom))$, and in order to estimate $\cn[\gga^{\ep}-\gga^{\eta};\cac_1^{\ga}]$, we shall resort to Lemma \ref{lem:4.3}  with $f=\gga^{\ep}-\gga^{\eta}$. This yields, for $p>1$ and $\ga<\al$,
$$
\cn[\gga^{\ep}-\gga^{\eta};\cac_1^{\ga}]\le c \, U_{\ga,p}\lp \delta(\gga^{\ep}-\gga^{\eta}) \rp
= c \lp \int_{\oom}\int_{\oom} \frac{|\delta(\gga^{\ep}-\gga^{\eta})_{ts}|^{2p}}{|t-s|^{2\ga p+4}}\, dsdt \rp^{1/2p}.
$$
Hence, invoking Jensen's inequality, we obtain:
\bean
\be\lc \cn[\gga^{\ep}-\gga^{\eta};\cac_1^{\ga}] \rc &\lesssim& 
 \lp  \int_{\oom}\int_{\oom} \frac{\be \lc  |\der(\gga^{\ep}-\gga^{\eta})_{ts}|^{2p}\rc}
{|t-s|^{2\ga p+4}} \, dsdt \rp^{1/2p}  \\
& \lesssim& 
 \lp  \int_{\oom}\int_{\oom} \frac{\be^{p} \lc  |\der(\gga^{\ep}-\gga^{\eta})_{ts}|^{2}\rc}
{|t-s|^{2\ga p+4}} \, dsdt \rp^{1/2p},
\eean
where we have used the fact that $\gga^{\ep},\gga^{\eta}$ are Gaussian processes in the last inequality. By considering $p$ large enough in the relation above, it is thus easily seen that, if we can prove that 
\begin{equation}\label{eq:63}
\be[  |\der(\gga^{\ep}-\gga^{\eta})_{ts}|^{2}]\le c |t-s|^{2\hat\ga} |\ep-\eta|^{\beta}
\end{equation} 
for a certain $\beta>0$ and $\ga<\hat\ga<\al$, 
then the following relation holds true:
$$
\be\lc  \cn[\gga^{\ep}-\gga^{\eta};\, \cac_1^{\ga}]\rc \lesssim |\ep-\eta|^{\beta}.
$$
Thus, we get that the family $\{\gga^{\ep};\,\ep>0\}$ is a Cauchy sequence in $L^1(\cu;\cac_1^\ga(\ott))$, whose limit is  the analytic fBm  $\gga$, {\em provided}  we can prove (\ref{eq:63}). The remainder of the proof is thus devoted to show the latter relation.

\smallskip

\noindent
{\it Step 2: Moment estimates.}
Set $U_{ts}=\be[  |\der(\gga^{\ep}-\gga^{\eta})_{ts}|^{2}]$. We shall now prove that $U_{ts}\le c_{\rho} 
|t-s|^{2\alpha\rho} |\ep-\eta|^{2\alpha(1-\rho)}$ for every $\rho\in(0,1)$. To this purpose, notice that:
$$
|\der(\gga^{\ep}-\gga^{\eta})_{ts}|^{2} \lesssim |\der \gga^{\ep}_{ts}|^{2} + |\der \gga^{\eta}_{ts}|^{2} ,
\text{ and }
|\der(\gga^{\ep}-\gga^{\eta})_{ts}|^{2} \lesssim |\gga^{\ep}_t-\gga^{\eta}_t|^{2} + |\gga^{\ep}_s-\gga^{\eta}_s|^{2} .
$$
This allows to write, for an arbitrary exponent $\rho\in(0,1)$,
$$
U_{ts} \lesssim \lp \be\lc |\der \gga^{\ep}_{ts}|^{2} + |\der \gga^{\eta}_{ts}|^{2} \rc \rp^{\rho}  \,
\lp \be\lc |\gga^{\ep}_t-\gga^{\eta}_t|^{2} + |\gga^{\ep}_s-\gga^{\eta}_s|^{2} \rc \rp^{1-\rho}.
$$
A direct application of Lemma \ref{lem:moments-B-ep} gives now:
\beq\label{eq:bnd-U-ts}
U_{ts} \lesssim |t-s|^{2\al\rho} \, |\ep-\eta|^{2\al(1-\rho)},
\eeq
which ends the proof, since $\ga\equiv\al\rho$ can be taken as close as we wish to $\al$.

\end{proof}

\begin{remark}\label{rmk:cvgce-B-ep-L-p}
A slight extension of the computations above allow to prove that in fact, $\gga^{\ep}$ converges in $L^p(\cu;\cac_1^\ga(\oom))$ for any $p>1$.
\end{remark}

\subsection{Convergence of Lévy's area}
\label{sec:levy-area}


Consider  a two-dimensional analytic fBm $\gga=(\gga(1),\gga(2))$ with independent components,
 and the associated approximation $\gga^{\ep}=(\gga^{\ep}(1),$ $\gga^{\ep}(2))$. We then  set 
\beq\label{eq:def-approx-levy}
\bde(j_1,j_2)=\int_{[s,t]} d\gga^{\ep}_{u_1}(j_1) \int_{[s,u_1]} d\gga^{\ep}_{u_2}(j_2),
\quad\mbox{for}\quad
j_1,j_2\in\{1,2\}, \, s,t\in\oom,
\eeq
where the above  iterated integral  is understood in the Riemann sense.
 This section is devoted to prove that $\bde$ is a convergent sequence 
in  $L^1(\cu;\cac_2^{2\ga}(\oom))$, and that its limit $\bd$ satisfies $\der \bd=\der B\otimes \der B$ as in Hypothesis \ref{hyp:x1}.
 we shall study the convergence for $j_1=j_2$ and $j_1\ne j_2$ separately.
\begin{proposition}\label{prop:cvgce-area-diag}
The increments $\bde(1,1)$ and $\bde(2,2)$ converge in $L^1(\cu;\cac_2^{2\ga}(\oom))$.
 \end{proposition}
 
 \begin{proof}
 For notational sake, we shall write $\cac_2^{2\ga},\cac_3^{2\ga}$ instead of $\cac_2^{2\ga}(\oom),\cac_3^{2\ga}(\oom)$ in the sequel. In order to prove that $\bde(1,1)$ is a Cauchy sequence in $L^1(\cu;\cac_2^{2\ga}(\oom))$, we invoke Proposition~\ref{prop:4.4}, which can be read here as:
\begin{multline*}
 \cn[\bde(1,1)-\bdet(1,1);\cac_2^{2\ga}]   \\
 \lesssim
 U_{2\ga,p}(\bde(1,1)-\bdet(1,1)) + \cn[\der\bde(1,1)-\der\bdet(1,1);\cac_3^{2\ga}]\equiv A_1+A_2.
 \end{multline*}
 In order to estimate the term $A_1$ above, notice first that, since $\gga^{\ep}$ is a regular path, we have: $\bde_{ts}(1,1)=\frac12 \lc \der \gga^{\ep}_{ts}(1) \rc^2$. Hence,
 $$
 A_1= \half \lp \int_{\oom}\int_{\oom} 
 \frac{|[\delta \gga^{\ep}_{ts}(1)]^2-[\delta \gga^{\eta}_{ts}(1)]^2|^{2p}}{|t-s|^{4\ga p+4}}\, dsdt \rp^{1/2p}.
 $$
 This integral can now be bounded as in Lemma \ref{lemma:cvgce-B-ep}, by means of an inequality similar to (\ref{eq:bnd-U-ts}).
 
 \smallskip
 
 Let us turn now to the evaluation of the term $A_2$. Owing to the fact that $\bde$ is defined by (\ref{eq:def-approx-levy}), where $\gga^{\ep}$ is a regular process, the following particular case of (\ref{eq:*}) is readily checked:
 $$
 \lc \der\bde(1,1)-\der\bdet(1,1)\rc_{tus}
 =\der \gga^{\ep}_{tu}(1) \, \der \gga^{\ep}_{us}(1) - \der \gga^{\eta}_{tu}(1) \, \der \gga^{\eta}_{us}(1).
 $$
 Hence, for an arbitrary coefficient $\rho\in(0,1)$ and $0<\ga<\hat\ga<\al$, we end up with:
 \begin{multline*}
  \lln \lc \der\bde(1,1)-\der\bdet(1,1)\rc_{tus} \rrn
  \le \lc  (\cn[\gga^{\ep};\cac_1^{\hat\ga}] + \cn[\gga^{\eta};\cac_1^{\hat\ga}])\, 
   |t-u|^{\hat\ga} |u-s|^{\hat\ga}\rc^{\rho} \\
  \times \Big[ \der \gga^{\ep}_{tu}(1) \lp \der \gga^{\ep}_{us}(1)-\der \gga^{\eta}_{us}(1) \rp 
  + \lp \der \gga^{\ep}_{tu}(1)-\der \gga^{\eta}_{tu}(1) \rp \, \der \gga^{\eta}_{us}(1)
 \Big]^{1-\rho},
 \end{multline*}
 and thus, a standard application of the Cauchy-Schwarz inequality yields:
 \begin{align*}
 \be  [A_2]&=\be\lc \cn[\der\bde(1,1)-\der\bdet(1,1);\cac_3^{2\ga}] \rc  \\
 &\lesssim \be^{1/2}\lc \cn^2[\gga^{\ep}(1);\cac_1^{\hat\ga}] + \cn^2[\gga^{\eta}(1);\cac_1^{\hat\ga}] \rc
 \, \be^{1/2}\lc \cn^{2(1-\rho)}[\gga^{\ep}(1)-\gga^{\eta}(1);\cac_1^{\hat\ga}] \rc  \\
 &\lesssim |\ep-\eta|^{\beta},
 \end{align*}
 for a certain $\beta>0$, according to Remark \ref{rmk:cvgce-B-ep-L-p}. Our claim is now easily deduced from our estimates on $A_1$ and $A_2$.
 
 \end{proof}
 
 \smallskip
 
 Let us begin the preliminary steps for the convergence of the crossed terms $\bde(1,2)$ and $\bde(2,1)$, for which the following notation will be needed:
\begin{notation}\label{not:F-ij}
For  $\ep_1,\ep_2>0$ and  $(x,\bar y)\in\bar\Pi^+\times\bar\Pi^-$, we set   
$$F_{\eps_1,\eps_2}(x,\bar y):=(-\imath(x-\bar y)+\eps_1+\eps_2)^{2\alpha-2}.
$$
\end{notation}
With this notation in hand, one can estimate the increments of $\bde$ as follows:
\begin{lemma}\label{lemma:bnd-incr-B2-ep}
For any $s,t\in\oom$ and $\rho\in(0,1)$, there exists a positive constant $c_{\rho}>0$ such that
$$
\be[ | (\bde(1,2)-\bdet(1,2))_{ts}|^{2}]\le c_{\rho} |t-s|^{2\alpha(1-\rho)} |\ep-\eta|^{2\alpha\rho}.
$$
\end{lemma}

\begin{proof}
According to identity (\ref{eq:cov-gamma-prime-path}), we have:
\begin{align*}
&\be[ | (\bde(1,2)-\bdet(1,2))_{ts}|^{2}]  \\
&=\be\Bigg[
\lp \ist  d\gga_{u_1+\im\ep}(1)  \int_{[s,u_1]} d\gga_{u_2+\im\ep}(2)
-\ist  d\gga_{u_1+\im\eta}(1)  \int_{[s,u_1]} d\gga_{u_2+\im\eta}(2) \rp\\
& \hspace{2.5cm}
\times \lp \ist  d\bar\gga_{v_1+\im\ep}(1)  \int_{[s,v_1]} d\bar\gga_{v_2+\im\ep}(2)
-\ist  d\bar\gga_{v_1+\im\eta}(1)  \int_{[s,v_1]} d\bar\gga_{v_2+\im\eta}(2) \rp
\Bigg] \\
&= \ist du_1 \int_{[s,u_1]} du_2 \ibst d\bar v_1 \int_{[\bar s,\bar v_1]} d\bar v_2 \,
F^{(2)}_{\ep,\eta}(u_1,\bar v_1;u_2,\bar v_2),
\end{align*}
where, recalling Notation \ref{not:F-ij}, the function $F^{(2)}_{\ep,\eta}$ is defined by:
\begin{align*}
&F^{(2)}_{\ep,\eta}(u_1,\bar v_1;u_2,\bar v_2)\\
&=
F_{\ep,\ep}(u_1,\bar v_1) F_{\ep,\ep}(u_2,\bar v_2)+ F_{\eta,\eta}(u_1,\bar v_1) F_{\eta,\eta}(u_2,\bar v_2)
-2  F_{\ep,\eta}(u_1,\bar v_1) F_{\ep,\eta}(u_2,\bar v_2) \\
&= F_{\ep,\ep}(u_1,\bar v_1) \lc F_{\ep,\ep}(u_2,\bar v_2)- F_{\ep,\eta}(u_2,\bar v_2) \rc
+ F_{\ep,\eta}(u_2,\bar v_2) \lc F_{\ep,\ep}(u_1,\bar v_1)-F_{\ep,\eta}(u_1,\bar v_1) \rc \\
& + F_{\eta,\eta}(u_1,\bar v_1) \lc F_{\eta,\eta}(u_2,\bar v_2)- F_{\eta,\ep}(u_2,\bar v_2) \rc
+ F_{\eta,\ep}(u_2,\bar v_2) \lc F_{\eta,\eta}(u_1,\bar v_1)-F_{\eta,\ep}(u_1,\bar v_1) \rc.
\end{align*}
We now have to control a sum made of many terms exhibiting the same level of difficulty. We shall thus focus on one of them, namely:
$$
I_1^{(2)}\triangleq
\ist du_1 \int_{[s,u_1]} du_2 \ibst d\bar v_1 \int_{[\bar s,\bar v_1]} d\bar v_2 \,
F_{\ep,\eta}(u_2,\bar v_2) \lc F_{\ep,\ep}(u_1,\bar v_1)-F_{\ep,\eta}(u_1,\bar v_1) \rc.
$$
For the control of $I_1^{(2)}$, as in Lemma \ref{lem:bound-intg-contour}, we introduce the contour of integration 
$$
\ga:=\ga_1\cup\ga_2\cup\ga_3=\big[ s,s+\imath|t-s|\big]\cup \big[s+\imath|t-s|,t+\imath|t-s| \big]
\cup \big[t+\imath|t-s|,t \big].
$$ 
If $z\in \gamma$, let $\gamma(z)$ be the section of the path $\gamma$ comprised between $s$ and $z$. Then
(by Cauchy's theorem)
$$
I_1^{(2)}=
\int_\ga dz_1 \int_{\ga(z_1)} dz_2 \int_{\bar{\ga}} d\bar{w}_1 \int_{\overline{\ga(w_1)}} d\bar{w}_2 \,
F_{\ep,\eta}(z_2,\bar{w}_2) \lc F_{\ep,\ep}(z_1,\bar{w}_1)-F_{\ep,\eta}(z_1,\bar{w}_1) \rc.
$$
As in Lemma \ref{lem:bound-intg-contour}, 9 terms should be controlled in order to achieve the desired bound. We shall treat the most divergent of them, that is:
$$
J_1^{(2)}=\int_{\ga_1} dz_1 \int_{\ga(z_1)} dz_2 \int_{\bar{\ga}_1} d\bar{w}_1 \int_{\overline{\ga(w_1)}}
 d\bar{w}_2 \,
F_{\ep,\eta}(z_2,\bar{w}_2) \lc F_{\ep,\ep}(z_1,\bar{w}_1)-F_{\ep,\eta}(z_1,\bar{w}_1) \rc.
$$
On $\ga_1$, the change of variable $z_1=s+\im u_1$ for $0\le u_1 \le |t-s|$, $z_2=s+\im u_2$ for $0\le u_2 \le u_1$, and the same kind of transformations for $\bar w_1,\bar w_2$, yield:
\begin{multline*}
J_1^{(2)}=
\int_0^{|t-s|} du_1 \int_0^{u_1} du_2 \int_0^{|t-s|} dv_1 \int_0^{v_1} d v_2  \\
\big[  (2\Im s+ u_1+v_1+2\ep)^{2\al-2}
-(2\Im s+u_1+v_1+\ep+\eta)^{2\al-2}\big] 
\, (2\Im s+u_2+v_2+2\ep)^{2\al-2},
\end{multline*}
and hence:
\begin{multline*}
|J_1^{(2)}|\le
\int_0^{|t-s|} du_1 \int_0^{|t-s|} dv_1 
\lln  (u_1+v_1+2\ep)^{2\al-2}-(u_1+v_1+\ep+\eta)^{2\al-2}\rrn  \\
\times \int_0^{|t-s|} du_2 \int_0^{|t-s|} dv_2 \lln (u_2+v_2+2\ep)^{2\al-2} \rrn.
\end{multline*}
As in Lemma \ref{lem:bound-intg-contour}, we can now easily conclude, for an arbitrary constant
 $0<\rho<1$, that:
$$
|J_1^{(2)}|\lesssim  |\ep-\eta|^{4\alpha\rho} (t-s)^{4\al(1-\rho)}.
$$
We may now  treat the other terms appearing in the analysis of $I_1^{(2)}$ (and more generally of $F^{(2)}$) in the same 
way, which ends the proof.

\end{proof}
 
 We are now ready to analyze the convergence of the crossed terms $\bde(1,2)$ and $\bde(2,1)$: 
 \begin{proposition}\label{prop:cvgce-area-cross}
The increments $\bde(1,2)$ and $\bde(2,1)$ converge in $L^1(\cu;\cac_2^{2\ga}(\oom))$.
 \end{proposition}
 
 \begin{proof}
 The beginning of the proof goes exactly along the same lines as for Proposition \ref{prop:cvgce-area-diag}.
Let us write $\cac_2^{2\ga},\cac_3^{2\ga}$ for $\cac_2^{2\ga}(\oom),\cac_3^{2\ga}(\oom)$. 
In order to prove that $\bde(1,2)$ is a Cauchy sequence in $L^1(\cu;\cac_2^{2\ga}(\oom))$, we invoke Proposition \ref{prop:4.4}:
\begin{multline*}
 \cn[\bde(1,2)-\bdet(1,2);\cac_2^{2\ga}]   \\
 \lesssim
 U_{2\ga,p}(\bde(1,2)-\bdet(1,2)) + \cn[\der\bde(1,2)-\der\bdet(1,2);\cac_3^{2\ga}]\equiv A_1+A_2.
 \end{multline*}
 The term $A_2$ can now be bounded as in Proposition \ref{prop:cvgce-area-diag}, owing to the fact that $\der\bde(1,2)=\der \gga^{\ep}(1)\, \der \gga^{\ep}(2)$. We thus get 
 $$
\be[A_2] \lesssim |\ep-\eta|^{\beta},
 $$
 for a certain $\beta>0$.
 
 \smallskip
 
The term $A_1$ can be handled in the following way: by definition, we have
 $$
 A_1= \lp \int_{\oom}\int_{\oom} 
 \frac{|\bde_{ts}(1,2)-\bdet_{ts}(1,2)|^{2p}}{|t-s|^{4\ga p+4}}\, dsdt \rp^{1/2p}.
 $$
We can now apply Jensen's inequality as in Lemma \ref{lemma:cvgce-B-ep}. Furthermore,  $\bde_{ts}(1,2)$ is a random variable in the second chaos of the fractional Brownian motion $\gga$, and since all the $L^p$ norms on any given fixed chaos are equivalent, we obtain:
 $$
 \be[A_1] \lesssim \int_{\oom}\int_{\oom} 
 \frac{\be^p\lc |\bde_{ts}(1,2)-\bdet_{ts}(1,2)|^{2}\rc}{|t-s|^{4\ga p+4}}\, dsdt.
 $$
 As in Lemma \ref{lemma:cvgce-B-ep}, we are now reduced to an estimate of the form 
 $$
 \be[  |(\bde-\bdet)_{ts}|^{2}]\le c |t-s|^{2\hat\ga} |\ep-\eta|^{\beta},
 $$ 
 for a certain $\beta>0$ and $\ga<\hat\ga<\al$. But this  estimate stems directly from Lemma~\ref{lemma:bnd-incr-B2-ep}, and gathering our estimates on $A$ and $B$, we have thus proved our claim:
 $$
 \be\lc \cn[\bde(1,2)-\bdet(1,2);\cac_2^{2\ga}] \rc
 \lesssim |\ep-\eta|^{\beta}.
 $$
 
 \end{proof}
 
 \begin{remark}
 As in Section \ref{sec:cvgce-B-ep}, the $L^p$-convergence of $\bde(i,j)$ in $\cac_2^{2\ga}(\Omega)$
 can also be obtained here for $i,j\in\{1,2\}$, by slightly adapting our computations for the $L^1$-convergence.
 \end{remark}
 
 Putting together the results we have obtained so far, we can now state the following existence result for a rough path of order 2 based on $\gga$, for any value of the Hurst parameter $\al\in(0,1/2)$:
 \begin{theorem} \label{th:convergence-2}
Let $\gga$ be an analytic fractional Brownian motion with Hurst parameter $\al\in(0,1/2)$, and $\gga^{\ep}$ its regular approximation.
 Let also $\bde$ be the regularized L\'evy area given by formula (\ref{eq:def-approx-levy}), and consider $0<\ga<\al$. Then:
 
\smallskip

\noindent
{\bf (1)}
For any $p\ge 1$, the couple $(\gga^{\ep},\bde)$ converges in $L^p(\cu;\cac_1^{\ga}(\oom;\R^d)\times\cac_2^{2\ga}(\oom;\R^{d,d}))$ to a couple $(\gga,\bd)$, where $\gga$ is the analytic fractional Brownian motion mentioned above.
 
\smallskip

\noindent
{\bf (2)}
The L\'evy area $\bd$ is an element of $\cac_2^{2\ga}(\oom;\R^{d,d})
\cap\cac_2^{{\rm m},\ga}(\oom;\R^{d,d})$.

\smallskip

\noindent
{\bf (3)}
The increment $\bd$ satisfies the multiplicative and geometric
algebraic relations prescribed in Hypothesis \ref{hyp:x}, namely: 
$$
\der\bd(i,j)=\der \gga(i) \, \der \gga(j),
\quad\mbox{and}\quad
\bd(i,j)+\bd(j,i)= \der \gga(i) \circ \der \gga(j),
$$
for $i,j=1,\ldots,d$.
 \end{theorem}
 
 \begin{proof}
 The first part of our assertion is trivially deduced from Propositions \ref{lemma:cvgce-B-ep}, \ref{prop:cvgce-area-diag} and \ref{prop:cvgce-area-cross}. The fact that $\bd$ is an element of $\cac_2^{{\rm m},\ga}(\oom;\R^{d,d}))$ can be shown thanks to a limiting procedure along the same lines as Propositions \ref{prop:cvgce-area-diag} and \ref{prop:cvgce-area-cross}, except that one has to replace the use of Proposition \ref{prop:4.4} by Lemma \ref{lem:4.3}.
 
 \smallskip
 
 As far as the third part of our claim is concerned, it is sufficient to notice that, since $\bde$ is a smooth process, the relation
 $$
\der\bde(i,j)=\der \gga^{\ep}(i) \, \der \gga^{\ep}(j),
\quad\mbox{and}\quad
\bde(i,j)+\bde(j,i)= \der \gga^{\ep}(i) \circ \der \gga^{\ep}(j),
$$
 is automatically satisfied, by some algebraic manipulations involving only usual Riemann integrals. The desired result is then obtained by taking limits on both sides of the identity above, and  taking into account that $\gga^{\ep}$ converges in any $L^q(\cu;\cac_1^{\ga}(\oom;\R^d))$, for $q\ge 1$.
 
 \end{proof}


\subsection{Multidimensional estimates}

\label{sec:multi-estim}

Let $\Gamma=(\Gamma(1),\ldots,\Gamma(d))$ be an $d$-dimensional analytic fractional Brownian motion.
This section is a generalization of the previous one to the case of multiply iterated integrals
of any order. As a result, we finally obtain a rough-path lying above $\Gamma$, which implies
the possibility to solve analytic stochastic differential equations driven by $\Gamma$ as in equation (\ref{eq:SDE}). Let us consider then our fixed bounded neighborhood $\oom$ of 0 and the analytic approximation $\gga^\ep$ of $\gga$. For $s,t\in\oom$, $n\le N=\lfloor 1/\al\rfloor$, and any tuple $(i_1,\ldots,i_n)\in\{1,\ldots,d\}^n$, the natural approximation of $\mathbf{\gga}^{\mathbf{n}}_{ts}(i_1,\ldots,i_n)$ is given by the Riemann iterated integral:
\begin{equation}
\mathbf{\gga}^{\mathbf{n},\ep}_{ts}(i_1,\ldots,i_n):=
\int_{[s,t]} d\gga^{\ep}_{u_1}(i_1) \int_{[s,u_1]} d\gga^{\ep}_{u_2}(i_2)\cdots 
\int_{[s,u_{n-1}]} d\gga^{\ep}_{u_n}(i_n). \label{eq:Gamma-n}
\end{equation}

In particular, we shall denote by ${\mathcal V}_{ts}^{\eps}$ the following L\'evy 'hypervolume':
$${\mathcal V}_{ts}^{\eps}=\mathbf{\gga}^{\mathbf{n},\ep}_{ts}(1,\ldots,n).$$

As in the case of Section \ref{sec:levy-area}, an important preliminary step in order to obtain the convergence of $\mathbf{\gga}^{\mathbf{n},\ep}$ is the following bound:
\begin{lemma}
Let $\Omega$ be a fixed bounded neighborhood of $0$ in $\bar{\Pi}^+$, and $p\ge 1$. For every
 $\rho\in(0,2\alpha)$, 
there exists a constant $C_{\rho}$ such that for every $n\ge 3$ and any $n$-uple $(i_1,\ldots,i_n)\in\{1,\ldots,d\}^n$, we have:
\begin{equation}\label{eq:67}
\be \lc\left| \mathbf{\gga}^{\mathbf{n},\ep}_{ts}(i_1,\ldots,i_n)
-\mathbf{\gga}^{\mathbf{n},\eta}_{ts}(i_1,\ldots,i_n)\right|^{2p}\rc
  \le C_{\rho} |\eps-\eta|^{p\rho} |t-s|^{p(2n\alpha-\rho)}, \quad s,t\in\Omega.
\end{equation}
\end{lemma}

\begin{proof}

First of all, since we are dealing with random variables in the $n\textsuperscript{th}$ chaos of the Gaussian process $\Gamma$, it is enough to prove inequality (\ref{eq:67}) for $p=1$. Next, the following lines prove that it is enough to estimate $\be [| \cv_{ts}^{\ep}-\cv_{ts}^{\eta}|^{2}]$.  Namely, suppose that some of the indices $(i_1,\ldots,i_n)$ coincide, and let $\Sigma_I\subset\Sigma_d$ be the subgroup of permutations $\sigma\in\Sigma_d$ such that
$i_{\sigma(j)}=i_j$ for all $j=1,\ldots,n$. Then
\begin{multline}\label{eq:70}
\be \lc\left| \mathbf{\gga}^{\mathbf{n},\ep}_{ts}(i_1,\ldots,i_n)
-\mathbf{\gga}^{\mathbf{n},\eta}_{ts}(i_1,\ldots,i_n)\right|^{2}\rc =\\
\sum_{\sigma\in \Sigma_I} \be \Big[ \big( \mathbf{\gga}^{\mathbf{n},\ep}_{ts}(1,\ldots,n)
-\mathbf{\gga}^{\mathbf{n},\eta}_{ts}(1,\ldots,n)\big) 
\overline{ \left(
\mathbf{\gga}^{\mathbf{n},\ep}_{ts}(\sigma(1),\ldots,\sigma(n))-
\mathbf{\gga}^{\mathbf{n},\eta}_{ts}(\sigma(1),\ldots,\sigma(n))
\right) } \Big], 
\end{multline}
and it is easily seen by the Cauchy-Schwarz inequality that this last term is bounded by $|\Sigma_I|\,\be [| \cv_{ts}^{\ep}-\cv_{ts}^{\eta}|^{2}]$. In order to justify  equation (\ref{eq:70}), let us just take the example $n=3$ and $(i_1,i_2,i_3)=(1,1,2)$. Then the computation of $\be [| \mathbf{\gga}^{\mathbf{3},\ep}_{ts}(1,1,2)
-\mathbf{\gga}^{\mathbf{3},\eta}_{ts}(1,1,2)|^{2}]$ involves products of the form:
$$
\be \lc\left(\Gamma'_{z_1}(1)\Gamma'_{z_2}(1)\Gamma'_{z_3}(2)\right)
\left(\bar{\Gamma}'_{w_1}(1)\bar{\Gamma}'_{w_2}(1)\bar{\Gamma}'_{w_3}(2)\right) \rc,
$$
which, invoking Wick's formula, are equal to
\begin{multline}\label{eq:71}
\be\lc \Gamma'_{z_1}(1)\bar{\Gamma}'_{w_1}(1)\rc \ 
\be \lc\Gamma'_{z_2}(1)\bar{\Gamma}'_{w_2}(1) \rc\ 
\be \lc\Gamma'_{z_3}(2)\bar{\Gamma}'_{w_3}(2) \rc \\
+\be \lc\Gamma'_{z_1}(1)\bar{\Gamma}'_{w_2}(1) \rc \ 
\be \lc\Gamma'_{z_2}(1)\bar{\Gamma}'_{w_1}(1) \rc\ 
\be \lc\Gamma'_{z_3}(2)\bar{\Gamma}'_{w_3}(2) \rc.
\end{multline}
These two terms correspond to  $\sigma=$Id or $\sigma=\tau_{12}$ in the right-hand side of (\ref{eq:70}), and one can check that expression (\ref{eq:71}) is equal to
$$
 \be \left[\left(\Gamma'_{z_1}(1)\Gamma'_{z_2}(2)\Gamma'_{z_3}(3)\right)
\left(\bar{\Gamma}'_{w_1}(1)\bar{\Gamma}'_{w_2}(2)\bar{\Gamma}'_{w_3}(3)
+ \bar{\Gamma}'_{w_1}(2)\bar{\Gamma}'_{w_2}(1)\bar{\Gamma}'_{w_3}(3)\right) \right].
$$
The general case can be treated along the same lines, up to some cumbersome notations.

\smallskip

Hence all we need is to  obtain a bound of the form:
\begin{equation}\label{eq:68}
\be \lc| \cv_{ts}^{\ep}-\cv_{ts}^{\eta}|^{2}\rc \le
C_{\rho} |\eps-\eta|^{\rho} |t-s|^{2n\alpha-\rho}, \quad s,t\in\Omega.
\end{equation}

\smallskip

As in the L\'evy area case of the previous section (see Lemma \ref{lemma:bnd-incr-B2-ep}), a straightforward application of identity (\ref{eq:cov-gamma-prime-path}) yields
\begin{multline*}
\be \lc| \cv_{ts}^{\ep}-\cv_{ts}^{\eta}|^{2}\rc 
= \left( \int_{[s,t]} dx_1 \int_{[s,x_1]} dx_2 \cdots \int_{[s,x_{n-1}]} dx_n\right) \\
\left( \int_{[\bar s,\bar t]} d\bar y_1 \int_{[\bar s,\bar y_1]} d\bar y_2 \cdots \int_{[\bar s,\bar x_{n-1}]} d\bar y_n\right) 
F^{(n)}_{\ep,\eta}(x_1,\bar y_1,\ldots,x_n,\bar y_n),
\end{multline*}
where, recalling Notation \ref{not:F-ij}, the function $F^{(d)}_{\ep,\eta}$ is defined by 
\begin{equation*}
F^{(n)}_{\ep,\eta}(x_1,\bar y_1,\ldots,x_n,\bar y_n)=
\prod_{j=1}^{n}F_{\ep,\ep}(x_j,\bar y_j) + \prod_{j=1}^{n}F_{\eta,\eta}(x_j,\bar y_j)
-2 \prod_{j=1}^{n}F_{\ep,\eta}(x_j,\bar y_j).
\end{equation*}
Observe that the latter function can be further decomposed into:
$$
F^{(n)}_{\ep,\eta}(x_1,\bar y_1,\ldots,x_n,\bar y_n)=
\sum_{j=1}^{d} G^j_{\ep,\eta}(x_1,\bar y_1,\ldots,x_n,\bar y_n)
+G^j_{\eta,\ep}(x_1,\bar y_1,\ldots,x_n,\bar y_n),
$$
where the functions $G^j$ are defined by:
\begin{multline*}
G^j_{\ep,\eta}(x_1,\bar y_1,\ldots,x_n,\bar y_n)=
F_{\ep,\ep}(x_1,\bar y_1)\cdots F_{\ep,\ep}(x_{j-1},\bar y_{j-1})
[F_{\ep,\ep}(x_j,\bar y_j)-F_{\ep,\eta}(x_j,\bar y_j)] \\
\times F_{\eps,\eta}(x_{j+1},\bar y_{j+1})\cdots F_{\eps,\eta}(x_n,\bar y_n)
\end{multline*}
We have thus proved that
$\be \lc| \cv_{ts}^{\ep}-\cv_{ts}^{\eta}|^{2}\rc
=\sum_{j=1}^d I^j_{\ep,\eta}+ I^j_{\eta,\ep}$,
where
\begin{multline*}
I^j_{\ep,\eta}=\left( \int_{[s,t]} dx_1 \int_{[s,x_1]} dx_2 \cdots \int_{[s,x_{n-1}]} dx_n\right) \\
\left( \int_{[\bar s,\bar t]} d\bar y_1 \int_{[\bar s,\bar y_1]} d\bar y_2 \cdots \int_{[\bar s,\bar x_{n-1}]} d\bar y_n\right) 
G^j_{\ep,\eta}(x_1,\bar y_1,\ldots,x_n,\bar y_n),
\end{multline*}
In order to show relation (\ref{eq:68}), it is thus sufficient to prove that, for all $j=1,\ldots,n$, we have 
$|I^j_{\ep,\eta}|\le C_{\rho} |\eps-\eta|^{\rho} |t-s|^{2n\alpha -\rho}$. Observe now that we may
 cast the term $I^j_{\ep,\eta}$ into the following form:
\begin{multline*}
I^j_{\ep,\eta}=\left( \int_{[s,t]} dx_1 \int_{[s,x_1]} dx_2 \cdots \int_{[s,x_{j-1}]} dx_j\right) 
\left( \int_{[\bar s,\bar t]} d\bar y_1 \int_{[\bar s,\bar y_1]} d\bar y_2 \cdots \int_{[\bar s,\bar x_{j-1}]} d\bar y_j\right) \\
F_{\eps,\eps}(x_1,\bar y_1)\cdots
F_{\eps,\eps}(x_{j-1},\bar y_{j-1}) \lc F_{\eps,\eps}(x_j,\bar y_j)-F_{\eps,\eta}(x_j,\bar y_j)\rc \phi(x_j,\bar y_j;s),
\end{multline*}
where
\begin{multline*}
 \phi(x_j,\bar y_j;s)=\left(\int_{[s,x_j]} dx_{j+1} \cdots 
\int_{[s,x_{n-1}]} dx_n\right)
\left(\int_{[\bar s,\bar y_j]} d\bar y_{j+1} \cdots \int_{[\bar s,\bar y_{n-1}]} d\bar y_n\right)\\
 F_{\eps,\eta}(x_{j+1},\bar y_{j+1})\cdots F_{\eps,\eta}(x_d,\bar y_d).
\end{multline*}
It is thus readily checked that the function $\phi(x_j,\bar{y}_j;s)$
 is an analytic iterated integral in the sense of Lemma \ref{lem:cvg-anal-intg},  bounded by
a constant times $|t-s|^{2\alpha(n-j)}$.
 Hence it
satisfies the hypothesis of  Lemma \ref{lem:bound-intg-contour}. As in the proof of the latter result, let 
$\gamma:=[s,s+\imath|t-s|]\cup[s+\imath|t-s|,t+\imath|t-s|]\cup[t+\imath|t-s|,t]$ be a complex deformation of the contour
 $[s,t]$, and, if $z\in \gamma$, let $\gamma(z)$ be the section of the path $\gamma$
comprised between $s$ and $z$. Then
\begin{multline*}
I^j_{\ep,\eta}=\left( 
\int_{\gamma} dx_1 \int_{\gamma(x_1)} dx_2 \ldots \int_{\gamma(x_{j-1})} dx_j 
\right) 
\left(
\int_{\bar\gamma} d\bar y_1 \int_{\overline{\gamma( y_1)}} d\bar y_2 \ldots \int_{\overline{\gamma( y_{j-1})}}
 d\bar y_j \right)\\
F_{\eps,\eps}(x_1,\bar y_1)\cdots
F_{\eps,\eps}(x_{j-1},\bar y_{j-1}) \lc F_{\eps,\eps}(x_j,\bar y_j)-F_{\eps,\eta}(x_j,\bar y_j)\rc \phi(x_j,\bar y_j;s),
\end{multline*}
and thus
\begin{multline*}
 |I^j_{\ep,\eta}| \le \left(
\int_{\gamma} |dx_1| \int_{\gamma} |dx_2| \ldots \int_{\gamma} |dx_j| \right)
\left(
\int_{\bar{\gamma}} |d\bar y_1| \int_{\bar{\gamma}} |d\bar y_2| \ldots \int_{\bar{\gamma}} |d\bar y_j| \right) \\
  |F_{\eps,\eps}(x_1,\bar y_1)|\cdots
|F_{\eps,\eps}(x_{j-1},\bar y_{j-1})| |F_{\eps,\eps}(x_j,\bar y_j)-F_{\eps,\eta}(x_j,\bar y_j)| |\phi(x_j,\bar y_j;s)|.
\end{multline*}
Now the multiple integral factorizes, and one is left with an expression of the form $A_1^{j-1}A_2$,
where

\begin{equation}
A_1=\int_{\gamma} |dx| \int_{\bar{\gamma}} |d\bar y| |(-\imath(x-\bar y)+2\eps)^{2\alpha-2}|\lesssim
|t-s|^{2\alpha}
\end{equation}

by Lemma \ref{lem:bound-intg-path}, and

\begin{equation}
A_2=\int_{\gamma} |dx_j| \int_{\bar{\gamma}} |d\bar{y}_j| |\phi(x_j,\bar y_j;s)|
\end{equation} 
which is bounded as in  Lemma \ref{lem:bound-intg-contour} by a constant times
$|t-s|^{2\alpha(n-j)+2\alpha-\rho} |\eps-\eta|^{\rho}$ for any $\rho\in(0,2\alpha)$.
 The above estimates now yield easily the desired bound
 $|I^j_{\ep,\eta}|\le C_{\rho} |\eps-\eta|^{\rho} |t-s|^{2n\alpha-\rho}$.

\end{proof}

The rough-path convergence of the multiplicative functional $(\Gamma^{\eps},\ldots,
\mathbf{\gga}^{\mathbf{n},\ep})$ to order $n$ is a consequence from the above computations and may
be stated as follows:

\begin{theorem}
Let $\gga$ be an analytic fractional Brownian motion with Hurst parameter $\al\in(0,1/2)$, and $\gga^{\ep}$ its regular approximation.
 Let also $n=\lfloor 1/\al\rfloor$, $\mathbf{\gga}^{\mathbf{k},\ep}$, $k=2,\ldots,n$ be the regularized iterated integrals
 given by formula (\ref{eq:Gamma-n}), and consider $0<\gamma<\al$. Then:
 
\smallskip

\noindent
{\bf (1)}
For any $p\ge 1$, the truncated multiplicative functional  $(\gga^{\ep},\mathbf{\gga}^{\mathbf{2},\ep},\ldots,
\mathbf{\gga}^{\mathbf{n},\ep})$
 converges in $L^p(\cu;\cac_1^{\ga}(\oom;\R^d)\times\ldots\times\cac_2^{n\ga}(\oom;(\R^{d})^{\otimes n}))$
 to an $n$-uple  $(\gga,\mathbf{\gga}^{\mathbf{2}},\ldots,\mathbf{\gga}^{\mathbf{n}})$,
 where $\gga$ is the analytic fractional Brownian motion defined in section 4.
 
\smallskip

\noindent
{\bf (2)}
The iterated integral $\mathbf{\gga}^{\mathbf{k}}$, $k=2,\ldots,n$
 is an element of $\cac_2^{k\gamma}(\oom;(\R^{d})^{\otimes k}))$, and it also belongs to $\cac_2^{{\rm m},\gamma}(\oom;
(\R^{d})^{\otimes k}))$.

\smallskip

\noindent
{\bf (3)}
The truncated multiplicative functional  $(\gga^{\ep},\mathbf{\gga}^{\mathbf{2},\ep},\ldots,
\mathbf{\gga}^{\mathbf{n},\ep})$ satisfies the multiplicative and geometric
algebraic relations prescribed in Hypothesis \ref{hyp:x}.
 \end{theorem}

\begin{proof}

Similar to that of Theorem \ref{th:convergence-2}.

\end{proof}


\end{document}